\documentclass{amsart}

\usepackage[T1]{fontenc}
\usepackage[british]{babel}
\usepackage{amssymb}
\usepackage{diagrams}

\makeatletter

\DeclareMathAlphabet\mathbfit{OML}{cmm}{b}{it}

\def\ie{i.\,e.}

\def\cf{cf.}

\def \N{{\bf N}}
\def \Z{{\bf Z}}

\def \C{{\bf C}}

\def\sign#1{\{#1\}}
\def\doublesign#1#2{\{#1,#2\}}
\def\degree#1{\lvert#1\rvert}

\DeclareMathOperator{\id}{id}

\DeclareMathOperator{\Hom}{Hom}

\def\dotcup{\mathbin{\dot\cup}}

\def \pair #1{{\langle #1\rangle}}

\def\t{\mathfrak t}

\def\M@warning#1{\typeout{MATH Warning: #1.}}

\def\M@margin#1{\if@printlabels\setbox0=\vbox to\z@
  {\vss\hbox to\z@{\hskip\hsize\hskip\labelskip\footnotesize #1\hss}}%
  \dp0=\z@\ifvmode\box0\else\vadjust{\box0}\fi\fi}

\def \M@restore{\catcode `\^=7 \catcode`\_=8 }

\def \newterm{\@ifnextchar[{\y@newterm}{\x@newterm}}
\def \x@newterm #1{\y@newterm[#1]{#1}}
\def \y@newterm[#1]#2{
  \textbf{\boldmath{#2}}}

\def \newsymbol #1{\M@margin{$\mapsto #1$}#1}
\def\printsymbols{{\def\indexname{Symbols}%
  \begin{theindex}%
  \normalsize
  \makeatletter
  \@input{\jobname.sym}%
  \makeatother
  \end{theindex}}}

\expandafter \ifx \csname diagram\endcsname \relax \else
\let \o@diagram=\diagram
\@namedef{diagram*}{\M@restore \o@diagram}
\@namedef{enddiagram*}{\enddiagram}
\def \diagram{\M@restore \refstepcounter{equation}%
        \o@diagram[eqno=\@eqnnum,moreoptions]}
\fi

\newif\if@printlabels
\let\printlabels=\@printlabelstrue
\let\printnolabels=\@printlabelsfalse
\newdimen\labelskip
\printnolabels
\let\o@label=\label
\newtoks\M@everypar
\def\M@label#1{\M@margin{\tt #1}\o@label{#1}%
  {\def\@currentlabel{\@currentthm}\o@label{lthm@#1}}%
  \@ifundefined{full@label}{}{\def\@currentlabel{\full@label}%
  \o@label{fthm@#1}}}
\def\label#1{\@ifundefined{@currentthm}{\o@label{#1}}%
    {\ifvmode\M@everypar=\everypar\everypar={\the\M@everypar\M@label{#1}}%
     \else\M@label{#1}\fi}}
\printnolabels

\let\o@enumerate\enumerate
\def\enumerate{\edef\full@label{%
  \@ifundefined{full@label}{\@currentlabel}{\full@label\,(\@currentlabel)}}%
  \o@enumerate}

\@ifpackageloaded{amsthm}{%
  \let\o@ynthm=\@ynthm  
  \def\@ynthm#1[#2]#3{\global\@namedef{dthm@#1}{#3}\o@ynthm{#1}[#2]{#3}}%
  }{%
  \let\o@othm=\@othm  
  \let\o@nthm=\@nthm
  \def\@othm#1[#2]#3{\global\@namedef{dthm@#1}{#3}\o@othm{#1}[{#2}]{#3}}%
  \def\@nthm#1#2{\global\@namedef{dthm@#1}{#2}\o@nthm{#1}{#2}}%
  }%
\let \o@thm=\@thm
\def \@thm{\let\@currentthm=\@currenvir\o@thm}
\def\@thmwarning#1{\M@warning{(no)theoremref: Reference `#1'
        on page \thepage \space not to a theorem}}
\def \@thmprintnum#1{\@ifundefined{r@fthm@#1}%
  {\ref{#1}}{\ref{fthm@#1}\,(\ref{#1})}}
\def\@thmcheck#1[#2]#3{\@ifundefined{r@#3}%
        {{\bf ??}\M@warning{(no)theoremref: Reference `#3' on page
        \thepage \space undefined}}%
        {\@ifundefined{r@lthm@#3}{\@thmwarning{#3}}%
        {\@ifundefined{dthm@#2}{\M@warning{(no)theoremref: Theorem
        environment `#2' on page \thepage \space undefined}}%
        {\edef \@tempa{\@nameuse{r@lthm@#3}}%
        \edef \@tempb{\expandafter \@car \@tempa \@nil}%
        \def \@tempa{#2}%
        \ifx \@tempa \@tempb \else
        \M@warning{(no)theoremref: Reference `#3' on page \thepage \space not
        to a `#2'}\fi}%
        #1\@thmprintnum{#3}\fi}}}
\def\@thmprintall#1{\@ifundefined{r@#1}%
        {{\bf ??}\M@warning{Reference `#1' on page \thepage \space undefined}}%
        {\x@thmprintall{#1}}}
\def \x@thmprintall#1{\@ifundefined{r@lthm@#1}%
        {\@thmwarning{#1}}{\y@thmprintall{#1}}}
\def \y@thmprintall#1{\edef \M@temp{\@nameuse{r@lthm@#1}}%
        \@nameuse{dthm@\expandafter \@car \M@temp \@nil}~\@thmprintnum{#1}}
\def\theoremref{\@ifnextchar[{\@thmcheck \iftrue}{\@thmprintall}}
\def\notheoremref{\@thmcheck{\iffalse}}

\def\smashsubstack#1{\hbox to5ex{\hss$\substack{#1}$\hss}}

\newenvironment{proofnoqed}[1][\proofname]{\par \normalfont
  \topsep6\p@\@plus6\p@ \trivlist \itemindent\normalparindent
  \item[\hskip\labelsep\scshape
    #1\@addpunct{.}]\ignorespaces
}{\endtrivlist}
\def\qedhere#1{\hbox to0pt{\hskip#1$\qed$\hss}}

\makeatother

\let\epsilon\varepsilon

\newtheorem{theorem}{Theorem}[section]
\newtheorem{lemma}[theorem]{Lemma}
\newtheorem{proposition}[theorem]{Proposition}
\newtheorem{corollary}[theorem]{Corollary}

\numberwithin{equation}{subsection}
\allowdisplaybreaks[1]

\newarrow{Equal}=====

\def\and{\qquad\text{and}\qquad}

\def\minusdegd{1}

\let\alg\mathbf
\let\Sim\mathbfit

\def\oneptspace{{\rm pt}}

\def\Ll{\boldsymbol\Lambda}
\def\Sl{\mathbf S}
\def\K{\mathbf K}
\def\Su{\Sl^*}
\def\Llstar{\Ll^{\!*}}
\def\Llprime{\Ll^{\!\prime}}

\def\p{\mathfrak p}
\def\pneg{\p\!}
\def\pzero{\mathfrak0\!}

\def\LMod{\hbox{\bf\em M\kern-0.05em od}}
\def\lMod#1{#1\hbox{-\bf\em M\kern-0.05em od}}
\def\RComod{\hbox{\bf\em Comod}}
\def\rComod#1{\hbox{{\bf\em Comod}-}#1}
\def\wMod#1{\lMod{B\Sl}}
\def\wComod#1{\rComod{B\Ll}}

\def\tpartial{\tilde\partial}
\let\shuffle\nabla
\def\AW{AW}
\def\AWtilde{\skew4\widetilde{\AW}}
\def\ST{ST}
\def\crossone{\mathbin{\times_1}}
\def\cupone{\mathbin{\cup_1}}
\def\crossoneproduct{cross$_1$~product}
\def\cuponeproduct{cup$_1$~product}
\def\Deltaone{\Delta^{(1)}}
\def\Bf{\tcomp_\oneptspace}
\def\pS{\tilde S}

\def\tcomp{\Psi}
\def\hcomp{\Phi}
\let\thom\psi
\let\hhom\phi

\def\T{\mathcal T}
\def\U{\mathcal U}

\def\LSpace{\hbox{\bf\em Action}}
\def\Space#1{#1\hbox{-{\bf\em Space}}}
\def\RSpaceover{\hbox{\bf\em Spaceover}}
\def\Spaceover#1{\hbox{{\bf\em Space}-}{#1}}

\def\Deltap{\tilde\Delta}
\def\Deltapp{\hat\Delta}

\def\timesunder#1{\mathbin{\mathchoice
  {\mathop\times\limits_{\mkern-5mu #1\mkern-5mu}}%
  {\times_{#1}}{\times_{#1}}{\times_{#1}}}}
\def\timesover#1{\mathbin{\mathchoice
  {\mathop\times\limits^{\mkern-3mu #1\mkern-5mu}}%
  {\times^{#1}}{\times^{#1}}{\times^{#1}}}}
\def\otimesunder #1{\mathbin{\mathchoice
  {\mathop\otimes\limits_{\mkern-20mu #1\mkern-20mu}}%
  {\otimes_{#1}}{\otimes_{#1}}{\otimes_{#1}}}}

\title[Koszul duality and equivariant cohomology for tori]%
  {Koszul duality\\and\\equivariant cohomology for tori}
\author{Matthias Franz}
\thanks{The author was partially supported by a grant
of the Deutscher Akademischer Austauschdienst.}

\address{Institut Fourier, Universit\'e de Grenoble~I, BP~74, 38402~Saint-Martin d'H\`eres, France}
\email{matthias.franz@uni-konstanz.de}


\begin{document}

\begin{abstract}
  Let $T$ be a torus.
  We show that Koszul duality can be used to compute the equivariant
  cohomology of topological $T$-spaces as well as the cohomology of pull backs
  of the universal $T$-bundle. The new features are
  that no further assumptions about the spaces are made and
  that the coefficient ring may be arbitrary.
  This gives in particular a Cartan-type model
  for the equivariant cohomology of a $T$-space
  with arbitrary coefficients.
  Our method works for intersection homology as well.
\end{abstract}

\maketitle

\tableofcontents


\section*{Introduction}

{ 

\def\theequation{\arabic{equation}}
\def\thetheorem{\ref{\linkedtheorem}${\mathbf '}\mkern-5mu$}

{ 
\def\alg#1{\mathop{\hbox{\em #1}}}
\def\t{\mathfrak t}

Let $T=(S^1)^r$ be a torus acting smoothly on some manifold~$X$.
The first definition of the $T$-equivariant cohomology~$H_T^*(X)$,
given by H.~Cartan, was as the homology of the differential graded
algebra
\begin{equation}\label{cartan-model}
  \Su\otimes\Omega^*(X)^T,
  \qquad
  d(\sigma\otimes\omega)=\sigma\otimes d\omega
    +\sum_{i=1}^r \xi_i\sigma\otimes x_i\cdot\omega,
\end{equation}
where $\Su$ is the algebra of polynomials on the Lie algebra~$\t$ of~$T$
and $\Omega^*(X)^T$ the complex of $T$-invariant differential forms on~$X$,
on which $\t$ acts by contraction with generating vector fields.
Moreover, $(x_i)$ denotes a basis of~$\t$
with dual basis~$(\xi_i)\in\Sl^2$.
(This so-called Cartan model exists, like all constructions we are about
to recall, for arbitrary compact connected Lie groups,
\cf~\cite{GuilleminSternberg:99}.
We are only concerned with the torus case, however,
and restriction thereto
will allow for a coherent presentation of previous results.)

Later on, A.~Borel defined equivariant cohomology for any topological
$T$-space and any coefficient ring~$R$
as the singular cohomology of the space~$ET\timesunder T X=(ET\times X)/T$,
where $ET$ is a contractible space on which $T$ acts freely.
For manifolds and real coefficients, this Borel construction gives
the same result as the Cartan model.

A few years ago, Goresky, Kottwitz
and MacPherson~\cite{GoreskyKottwitzMacPherson:98}
related equivariant cohomology to Koszul duality
as defined by Bernstein, Gelfand and Gelfand~\cite{BernsteinGelfandGelfand:78}:
They remarked that the complex~\eqref{cartan-model} arises from~$\Omega^*(X)^T$
by applying the Koszul functor~$\alg t$,
which carries differential graded modules over~$\Ll_*=\bigwedge^*\t=H_*(T)$,
the homology of the torus, to those over~$\Su=H^*(BT)$,
the cohomology of the classifying space of~$T$.
(The multiplication of~$\Ll_*$ is induced by that of~$T$.)
They went on to show that application of the quasi-inverse
Koszul functor~$\alg h$ to a suitably defined
differential $\Su$-module
of differential forms on~$ET\timesunder T X$ gives
a differential $\Ll_*$-module
\begin{equation}\label{complex-h-omega}
  \Omega^*(ET\timesunder T X)\otimes\Llstar,
  \qquad
  d(\omega\otimes\alpha)=d\omega\otimes\alpha
    +(-1)^{\degree\omega}\sum_{i=1}^r\,\xi_i\cdot\omega\otimes x_i\cdot\alpha
\end{equation}
whose homology is isomorphic to~$H^*(X)$ as $\Ll_*$-module.
(Here $\Llstar=H^*(T)$ is the cohomology of~$T$.)
They actually prove a much more general statement that for instance
also subsumes intersection homology of subanalytic spaces.

Shortly afterwards, Allday and Puppe~\cite{AlldayPuppe:99} pointed out
that the appearance of Koszul duality here reflects
an underlying topological duality between
the category~$\Space T$ of $T$-spaces on the one hand
and the category~$\Spaceover{BT}$ of spaces over~$BT$
on the other. (A space~$Y$ over~$BT$ is map~$Y\to BT$.)
The Borel construction is a functor
\[
  \Sim t\colon\Space T\to\Spaceover{BT},
\]
and pulling back the universal $T$-bundle~$ET\to BT$ along~$Y\to BT$
gives a functor in the other direction,
\[
  \Sim h\colon\Spaceover{BT}\to\Space T.
\]
From this point of view Goresky, Kottwitz and MacPherson's result for
ordinary cohomology can be stated as follows:
\begin{subequations}\label{G-K-M-result-A-P}
\begin{align}
  H^*\bigl(\alg t\Omega^*(X)^T\bigr) &\cong H^*\bigl(\Omega^*(\Sim t X)\bigl)
  & \hbox{as $\Su$-modules}\\
  H^*\bigl(\alg h\Omega^*(Y)\bigr) &\cong H^*\bigl(\Omega^*(\Sim h Y)^T\bigr)
  & \hbox{as $\Ll_*$-modules}
\end{align}
\end{subequations}
for smooth~$X\in\Space T$~and~%
~$Y=\Sim t X$.
(We have~$H^*(\Sim h Y)=H^*(\Sim h \Sim t X)\cong H^*(X)$
because $\Sim h\Sim t X\cong ET\times X$
is homotopy-equivalent to~$X$.)

The main purpose of this paper is to generalise this to an arbitrary
coefficient ring~$R$ instead of the reals.
We will use singular homology (and cohomology) rather than subanalytic chains,
and this will allow us to drop the assumption of subanalyticity
made in~\cite{GoreskyKottwitzMacPherson:98}.

Note that simply replacing de~Rham complexes by singular cochains
in~\eqref{G-K-M-result-A-P} does not make sense: While it is easy to define
a $\Ll_*$-action on~$C^*(X;R)$, $X$ a $T$-space, it is not possible in general
to lift the $\Su$-action on the cohomology of a space~$Y$ over~$BT$
to~$C^*(Y;R)$.
The reason for this is that the cup product of singular cochains
is inherently non-commutative.
(It is only commutative up to homotopy, hence commutative
in cohomology.) As a consequence, if one pulls back representatives~$\xi_i'$
of the generators~$\xi_i\in\Su=H^*(BT;R)$ along~$p\colon Y\to BT$
and imitates definition~\eqref{complex-h-omega}
with~$\xi_i\cdot\gamma=\gamma\cup p^*(\xi_i')$ for~$\gamma\in C^*(Y;R)$,
then $d$ will be not a differential any more.
(The only reason to multiply by~$p^*(\xi_i')$ from the \emph{right}
is to be consistent with other constructions later on.)

Since we ultimately only want \eqref{complex-h-omega} to be a differential
$\Ll_*$-module, the key idea is to perturb the differential
by ``higher order terms'':
\begin{equation}\label{gugenheim-may-differential-introduction}
  d(\gamma\otimes\alpha)=d\gamma\otimes\alpha
    -\sum_{\smashsubstack{\pi\subset\{1,\ldots,r\}\\\pi\ne\emptyset}}
      (-1)^{\degree\gamma+\degree\pi}\,
      \gamma\cup p^*(\xi_\pi')\otimes x_\pi\cdot\alpha.
\end{equation}
Here $\xi_\pi'\in C^{\degree\pi+1}(BT;R)$, and
$(x_\pi)$ denotes the canonical $R$-basis of~$\Ll_*$ generated by the~$x_i$'s.
This formula already appears in the work
of Gugenheim and May~\cite{GugenheimMay:74}, and they show how to
construct some~$(\xi_\pi')$ such that the resulting map~$d$
is a differential (which then is necessarily $\Ll_*$-equivariant).
They also prove that the complex~$C^*(Y;R)\otimes\Llstar$ with this new
differential indeed computes the cohomology of the pull-back~$\Sim h Y$,
but their method only gives an isomorphism of $R$-modules, \ie, without the
additional $\Ll_*$-action.

Choosing a collection~$(\xi_\pi')$ as above actually means defining a so-called
twisting cochain~$t\in\Hom_{-1}(C_*(BT;R),\Ll_*)$, which allows one to define
twisted tensor products like~\eqref{cartan-model},~\eqref{complex-h-omega}
or~\eqref{gugenheim-may-differential-introduction} in a conceptual manner.
It is also tantamount to a morphism of differential coalgebras from~$C_*(BT;R)$
to~$B\Ll_*$, the reduced bar construction of~$\Ll_*$. Comodules over~$B\Ll_*$
are also called $\Su$-modules `up to homotopy', and this hints to the fact
(well-known to experts)
that one can extend Koszul duality to these objects.

} 

\medbreak

The results of this paper are most naturally formulated in the homological
setting; this also facilitates the passage to intersection homology later on.
So we look at $C_*(Y;R)$ as an $\Sl_*$-comodule up to homotopy, which we call
a `weak $\Sl_*$-comodule'. Here $\Sl_*=H_*(BT;R)$ is the homology coalgebra of~$BT$.
This implies that the Koszul functors $\alg t$~and~$\alg h$
of this paper are \emph{not}
the ones from~\cite{GoreskyKottwitzMacPherson:98}, which we have used above,
nor those from~\cite{BernsteinGelfandGelfand:78}.
The functor~$\alg h$ carries weak $\Sl_*$-comodules to $\Ll_*$-modules,
and in order to be symmetric, we extend $\alg t$ to analogously defined
weak $\Ll_*$-modules, which are transformed into $\Sl_*$-comodules.
Any (strict) $\Ll_*$-module or $\Sl_*$-comodule has a canonical weak (co)module
structure, and we consider the Koszul functors as mapping
to the categories of weak (co)modules.

The main novelty of the present paper is the construction of
natural transformations
\begin{subequations}\label{connecting-morphisms}
\begin{align}
  \tcomp\colon\alg t\circ C_* &\to C_*\circ\Sim t\\
  \hcomp\colon C_*\circ\Sim h &\to \alg h\circ C_*.
\end{align}
\end{subequations}
Here we have used the same symbol~$C_*$ to denote both the functor assigning
the $\Ll_*$-module~$C_*(X;R)$ to a $T$-space~$X$
and the functor assigning to a space~$Y$ over~$BT$
the weak $\Sl_*$-comodule~$C_*(Y;R)$.

Once these natural transformations are established, an easy application
of the Leray--Serre theorem will prove our main result, which generalises
\eqref{G-K-M-result-A-P} to arbitrary topological spaces and arbitrary
coefficient ring:

\def\linkedtheorem{natural-transformations-c-equivalences}
\begin{theorem}\label{main-theorem-introduction}
  The natural transformations~$\tcomp$~and~$\hcomp$ are quasi-equivalences,
  \ie, they induce isomorphisms in homology for all objects.
\end{theorem}

It turns out that in the case of the Borel construction one can do better.
The cohomological formulation is a complete generalisation
of the Cartan model to singular cochains:

\def\linkedtheorem{singular-Cartan-model}
\begin{theorem}\label{singular-Cartan-model-introduction}
  Let $\Su$ act on~$\Su\otimes C^*(X;R)$ in the canonical way
  and give the tensor product
  the following differential and $\Su$-bilinear product:
  \begin{align*}
    d(\sigma\otimes\gamma) & =\sigma\otimes d\gamma
      +\sum_{i=1}^r\,\xi_i\sigma\otimes x_i\cdot\gamma, \\
   (\sigma'\otimes\gamma')(\sigma\otimes\gamma)
      &= \sigma'\sigma\otimes\gamma'\cup\gamma.
  \end{align*}
  Then $\tcomp_X^*\colon C^*(ET\timesunder T X)\to\Su\otimes C^*(X)$
  is a morphism of algebras and induces an $\Su$-equivariant isomorphism
  in homology.
\end{theorem}

Suppose that $X$ is a stratified pseudomanifold and $\p$~a perversity.
In this case we define a certain $\Ll_*$-submodule~$N_*$ of the singular
chain complex~$C_*(X;R)$ and show in a further step:

\def\linkedtheorem{allowable-subsets-theorem}
\begin{theorem}\label{main-theorem-2-introduction}
  The morphism~$\tcomp_X$ induces an
  isomorphism~$H(\alg t N_*)\cong I^{\pneg}H^T_*(X;R)$.
  A similar statement holds for
  spaces over~$BT$.
\end{theorem}

(See the main text for a precise statement of each theorem.)

\medbreak

Equivariant cohomology has attracted much interest from other parts of
mathematics in recent years, for example from symplectic and algebraic geometry
or combinatorics.
I have therefore striven to keep this paper reasonably self-contained
and accessible to non-experts.
It is organised as follows:

In Chapter~\ref{algebra} we introduce weak $\Sl_*$-comodules and
weak $\Ll_*$-modules and sketch the proof of
Koszul duality between these categories. This actually holds in much greater
generality and is usually formulated in the framework of operads,
see for example~\cite[Sec.~4.2]{GinzburgKapranov:94},~%
\cite[\S\,5.1]{Loday:96}, or~\cite{KrizMay:95}.
We nevertheless present an elementary version in order to save the reader from
digesting a sophisticated formalism before arriving at the relatively easy
result needed here.

Instead of using topological constructions of $ET$~and~$BT$, we will work in
the simplicial category \cite{May:68},~\cite{Lamotke:68}.
In the last section of the second chapter we define
the topological counterparts of the Koszul functors
in the simplicial setting and show that they are quasi-inverse to each other
-- not only for tori, but for arbitrary topological or even simplicial groups.
I propose to call these functors ``simplicial Koszul functors''
in order to stress the similarity
between them and the usual ``algebraic'' Koszul functors.
The preceding Sections \ref{notation-2}~to~\ref{classifying-spaces}
are essentially a review of well-known results about simplicial sets,
and except for Propositions
\theoremref[proposition]{properties-Steenrod}~and~%
\theoremref[proposition]{cone-homotopy}, I make no claim of originality.
One reason for including the material here is again to make the whole paper
readable for the non-specialist.
More importantly, we shall make essential use of the exact form
of each definition given in this chapter.
In most cases one can find several slightly different versions
in the literature, and an  underlying ``theorem'' of this article is
that they can all be chosen consistently.

Chapter~\ref{main-chapter} is the heart of the paper.
After explaining how to view the chain complex of a $T$-space
as a $\Ll_*$-module (which is easy) and the chain complex of a space~$Y$
over~$BT$ as a weak $\Sl_*$-comodule (Gugenheim--May),
we define the natural transformations~\eqref{connecting-morphisms}
and prove the main theorem.

The singular Cartan model is presented in Chapter~\ref{cohomology}.
For~$X=\oneptspace$~a point our theory gives a quasi-isomorphism
of algebras~$C^*(BT;R)\to\Su=H^*(BT;R)$.
We prove that this map has the important property of annihilating
all \cuponeproduct s. Our construction of a such map is much shorter
than the original one given in~\cite{GugenheimMay:74}.

In Chapter~\ref{intersection-homology} we generalise
\theoremref{natural-transformations-c-equivalences}
to intersection homology. Starting with the geometric definition of
intersection homology given in~\cite{GoreskyMacPherson:80}, we use
results due to Goresky--MacPherson~\cite{GoreskyMacPherson:86}
and  King~\cite{King:85} to arrive at the notion of an ``allowable subset''
of a simplicial set, which is defined as a graded subset closed
under all degeneracy maps and under the last face map.
Inspection of the proof of the main theorem then shows that it holds true
for allowable subsets, hence for intersection homology.

\medbreak

The core of the paper is elementary. Rudimentary knowledge of
homological algebra and simplicial sets is amply sufficient,
granted the occasional use of spectral sequences, in particular
of the Leray--Serre spectral sequence of a fibre bundle.
Most proofs are straightforward verifications
that some claimed identities do hold and are often only sketched
or entirely left to the reader.

Apart from reading this paper from beginning to end,
there are at least two more meaningful ways:
Readers only interested in the singular Cartan model should check
Sections~\ref{notation}~and~\ref{notation-2} for notation and should then,
after a quick look at Section~\ref{products},
turn directly to Section~\ref{section-Cartan}.
Those who would like to learn more about the relationship
between algebraic and simplicial Koszul duality, but without
entering into the details of how to construct the
natural transformations~\eqref{connecting-morphisms}
may skip Sections
\ref{section-Steenrod},~\ref{fibre-bundles},~\ref{classifying-spaces}
and \ref{important-maps-1},~\ref{important-maps-2}.

} 

\medbreak

\emph{Acknowledgements.}\hskip1em
This paper grew out of my thesis~\cite{Franz:01},
and I heartily thank my advisor
Volker Puppe for suggesting this topic to me and also for many inspiring
discussions.
Substantial additions were made during a visit at MIT in spring~2002,
and I thank the MIT~maths department
and in particular Victor Guillemin for their hospitality.
I am also indebted to Mark Goresky for a discussion of the sweep action,
to Robert MacPherson for pointing out King's definition
of intersection homology to me,
to the referee for many valuable comments,
and to Ezra Miller, Tara Holm and Ludger Kaup
for stylistic improvements.

\section{Algebraic Koszul duality}\label{algebra}


\subsection{Notation and terminology}\label{notation}

Throughout this article the letter~$R$ denotes
a commutative ring with unit element.

The set~$\{1,\ldots,r\}$ is denoted by~$[r]$.
The degree of a finite set is the number of its elements.
Generally, whenever an object~$x$ has a degree~$\degree x\in\Z$,
we call $\sign x=(-1)^{\degree x}$ its \newterm{sign};
moreover, we write $\doublesign x y=(-1)^{\degree x\cdot\degree y}$.
A summation over~$(\mu,\nu)\vdash\pi$ extends over all
partitions~$\mu\dotcup\nu=\pi$, and one over~$(\mu,\nu)\vdash(m,n)$
over all $(m,n)$-shuffles, \ie, over all
partitions~$\mu\dotcup\nu=\{0,\ldots,m+n-1\}$
with $\degree\mu=m$~and~$\degree\nu=n$.
We write $\sign{(\mu,\nu)}$ for the sign of the permutation determined by
such a partition.

The following definitions and sign conventions for complexes
seem to be consistent with~\cite{Dold:80} and the nice summary in
the first sections of~\cite[Ch.~3]{Smirnov:01}:

If not otherwise specified, a \newterm{complex} is one of $R$-modules,
and all tensor products and homomorphism complexes are taken over~$R$.
The $n$-th degree of a complex~$N$ is denoted by~$N_n$.
(Note that from now on we suppress meaningless dots as in~``$N_*$''.)
Let $M$,~$N$ be complexes. The complex~$M\otimes N$ has
differential~$d(m\otimes n)=d m\otimes n+\sign m\,m\otimes d n$,
and that on~$\Hom(M,N)$ is $d(f)(m)=d f(m)-\sign f\,f(d m)$.
(The former differential illustrates the general sign rule to insert,
whenever two objects are transposed, the factor $-1$ to the product
of their degrees. The latter is the \emph{only} exception.)
The canonical chain map
\begin{equation}\label{isomorphism-Hom-tensor}
  N\otimes M^*\to\Hom(M,N),
  \qquad
  a\otimes\gamma\mapsto\Bigl(c\mapsto\pair{\gamma, c}\,a\Bigr),
\end{equation}
is an isomorphism if $M$~or~$N$ is a finitely generated free $R$-module.
We denote by~$T_{MN}$ the canonical
transposition~$M\otimes N\to N\otimes M$,
$m\otimes n\mapsto\doublesign m n\,n\otimes m$.
Moreover, if $f\colon M\to M'$~and~$g\colon N\to N'$ are maps of complexes,
then $f\otimes g\colon M\otimes N\to M'\otimes N'$ is defined
by~$(f\otimes g)(m\otimes n)=\doublesign g m\,f(m)\otimes g(n)$.

By an \newterm{algebra}, we mean an associative differential graded
$R$-algebra~$A$ with multiplication~$\mu_A\colon A\otimes A\to A$
and unit~$\iota_A\colon R\to A$, and by a \newterm{coalgebra}, a coassociative
differential graded $R$-coalgebra~$C$ with
comultiplication~$\Delta_C\colon C\to C\otimes C$ and
augmentation~$\epsilon_A\colon C\to R$.
We analogously  write $\mu_N$~and~$\Delta_M$ for the structure maps of a
(left or right) $A$-module~$N$ and a (left or right) $C$-comodule~$M$,
respectively.
We call $A$ commutative if~$\mu_A T_{AA}=\mu_A$, and $C$ cocommutative
if~$T_{CC}\Delta_C=\Delta_C$.
(This is often called ``graded (co)commutative''.)
The category of left $A$-modules is denoted by~$\lMod A$, and that of right
$C$-comodules by~$\rComod C$.

Let $C$ be a coalgebra and $A$ an algebra.
Then $\Hom(C,A)$ is an algebra with \newterm{cup product}
\begin{equation}\label{definition-cup-product}
  u\cup v=\mu_A(u\otimes v)\Delta_C
\end{equation}
and unit~$\iota_A\epsilon_C$.
This applies in particular to~$C^*=\Hom(C,R)$.
Note that a map~$\gamma\colon C_n\to R$ has degree~$-n$, and that
the differential on~$C^*$, like all differentials,
\emph{lowers} degree, \cf~\cite[Ex.~1.9]{Dold:80}.
In order not to confuse the reader too much, we introduce the
notation~$C^n=(C^*)_{-n}$, so that $\gamma\in C^n$.

A right $C$-comodule~$M$ is canonically a left $C^*$-module by the
\newterm{cap product}
\begin{equation}\label{definition-cap-product}
  \gamma\cap m=(1\otimes\gamma)\Delta_M(m)\in M\otimes R=M
\end{equation}
for~$\gamma\in C^*$~and~$m\in M$.

We shall occasionally deal with coalgebras whose structure maps are
equivariant with respect to some Hopf algebra~$A$.
(Recall that the tensor product~$N\otimes N'$ of two $A$-modules
is again an $A$-module.)
We call such a coalgebra an \newterm{$A$-coalgebra}.

\subsection{Twisting cochains}\label{twisting}

Let $A$ be an algebra and $C$ a coalgebra.
A \newterm{twisting cochain}~\cite{Brown:59} is
an element~$u\in\Hom_{-1}(C,A)$ such that $d(u)+u\cup u=0$.
(See for example also~\cite[Sec.~3.3]{Smirnov:01} for the following.)

If $M$ is a right $C$-comodule and $N$ a left $A$-module, then the
\newterm{twisted tensor product}~$M\otimes_u N$ has the same underlying
graded $R$-module as the usual tensor product~$M\otimes N$,
but with differential
\begin{equation}
  d_u=d_M\otimes1+1\otimes d_N
      +(1\otimes\mu_N)(1\otimes u\otimes1)(\Delta_M\otimes1).
\end{equation}
Note that $C\otimes_u N$ is canonically a left $C$-comodule
and $M\otimes_u A$ a right $A$-module.

Similarly, for a left $C$-comodule~$M$ and a left $A$-module~$N$
we define the \newterm{twisted complex of homomorphisms}~$\Hom^u(M,N)$
with differential
\begin{equation}
  d^u(s)=d(s)+\mu_N(u\otimes s)\Delta_M.
\end{equation}

\subsection{The Koszul complex}\label{Koszul-complex}

Let $P$ be a free graded $R$-module of finite rank~$r$ which is concentrated
in positive odd degrees.
(If the characteristic of~$R$ is~$2$, then signs do
not matter, and $P$ may also have elements in even positive degrees.)
Denote by $\Ll=\bigwedge P$ the exterior algebra over~$P$ and
by~$\Sl$ the symmetric coalgebra over~$P[-1]$.
($P[-1]$ is obtained from~$P$ by raising all degrees by~$1$.)
We shall use that $\Ll$ is a Hopf algebra with augmentation~$\epsilon_{\Ll}=0$
and comultiplication~$\Delta_{\Ll}(x)=x\otimes1+1\otimes x$
for all~$x\in P\subset\Ll$.
For~$\Sl$ we need the canonical unit~$\iota_{\Sl}\colon R\to\Sl$.
The coefficient ring~$R$ is a left $\Ll$-module
via~$\epsilon_{\Ll}\colon\Ll\otimes R\to R$
and a right $\Sl$-comodule via~$\iota_{\Sl}\colon R\to R\otimes\Sl$.

Choose a basis $(x_i)$ of~$P$ with dual basis~$(\xi_i)$ of~$P^*$.
This gives canonical bases $(x_\pi)$ of~$\Ll$
and $(\xi^\alpha)$ of the algebra~$\Su$ dual to~$\Sl$ with
\[
  x_\pi=x_{i_1}\wedge\cdots\wedge x_{i_q}\in\Ll,
  \qquad
  \xi^\alpha=\xi_1^{\alpha_1}\cdots \xi_r^{\alpha_r}\in\Su,
\]
where $\pi=\{i_1<\cdots<i_q\}\subset[r]$ and~$\alpha\in\N^r$.
(Take this as the definition of~$\Ll$ if the characteristic of~$R$ is~$2$.)
In terms of the $R$-basis~$x^\alpha\in\Sl$ dual to the~$\xi^\alpha$'s,
the comultiplication of~$\Sl$ takes the form
\begin{equation}
  \Delta_{\Sl}(x^\alpha)=\sum_{\beta+\gamma=\alpha}x^\beta\otimes x^\gamma.  
\end{equation}

There is a canonical twisting cochain
\begin{equation}\label{definition-canonical-twisting-cochain}
  u_P = \sum_{i=1}^r x_i\otimes\xi_i\in\Hom_{-1}(\Sl,\Ll),
\end{equation}
which does not depend on the chosen bases.
(Here we have used the isomorphism~\eqref{isomorphism-Hom-tensor}.)
We write $M\otimes_P N$ for the tensor product with twisting cochain~$u_P$.

The \newterm{Koszul complex}
\begin{equation}\label{Koszul-differential}
  \newsymbol{\K}=\Sl\otimes_P\Ll,
  \quad
  d(s\otimes a)=\sum_{i=1}^r\xi_i\cap s\otimes x_i\wedge a
\end{equation}
is a left $\Sl$-comodule and a right $\Ll$-module. It is actually
a $\Ll$-coalgebra with componentwise comultiplication
\[
  \Delta_{\K}=(1\otimes T_{\Sl\Ll}\otimes1)(\Delta_{\Sl}\otimes\Delta_{\Ll})
  \colon
  \K=\Sl\otimes_P\Ll\to(\Sl\otimes_P\Ll)\otimes(\Sl\otimes_P\Ll)=\K\otimes\K
\]
and augmentation~%
$\epsilon_{\K}=\epsilon_{\Sl}\otimes\epsilon_{\Ll}\colon\K\to R$.

The latter map 
as well as the morphism of $\Sl$-modules
$\iota_{\Ll}\otimes\iota_{\Sl}\colon R\to\K$
are quasi-isomorphisms, \ie, they induce isomorphisms in homology.

\subsection{Weak $\Sl$-comodules}\label{weak-comodules}

We define a \newterm{weak \textmd{(right)} $\Sl$-comodule} to be a triple
$(C, M, u)$,
where $C$ is a coalgebra, $M$ a right $C$-comodule,
and $u\in\Hom_{-1}(C,\Ll)$ a twisting cochain with~$\epsilon_{\Ll}u=0$.

Using again the isomorphism~\eqref{isomorphism-Hom-tensor},
we may express $u$ in the form
\begin{equation}
    u=\sum_{\emptyset\ne\pi\subset[r]}x_\pi\otimes\gamma_\pi
\end{equation}
for some~$\gamma_\pi\in C^{\degree\pi+1}=(C^*)_{-\degree\pi-1}$.
The condition of $u$ being a twisting cochain translates into the equations
\begin{equation}\label{equation-a-pi}
  \forall\,\emptyset\ne\pi\subset[r]\qquad
    d\gamma_\pi=-\sum_{\smashsubstack{(\mu,\nu)\vdash\pi \\
                                      \mu\ne\emptyset\ne\nu}}
    \sign\mu\sign{(\nu,\mu)}\,\gamma_\mu\cup\gamma_\nu,
\end{equation}
and the explicit form of the differential on~$M\otimes_u\Ll$ is
\begin{equation}\label{differential-M-Lambda}
  d(m\otimes a)=d m\otimes a
    +\sum_{\pi\ne\emptyset}\doublesign\pi m\,
       \gamma_\pi\cap m\otimes x_\pi\wedge a.
\end{equation}

For small~$\pi$, condition~\eqref{equation-a-pi} says
$d\gamma_i = 0$, and~$d\gamma_{i j} = \gamma_j\cup\gamma_i-\gamma_i\cup\gamma_j$
for~$i<j$. This allows us to define a right $\Sl$-comodule structure on~$H(M)$,
or equivalently, a left $\Su$-module structure,
by setting
\begin{equation}\label{definition-module-structure-homology}
  \xi_i\cdot[m]=[\gamma_i\cap m].
\end{equation}
(This is equivalent because $\Sl^{**}=\Sl$.)
Hence while $M$ is only an $\Sl$-comodule `up to homotopy',
$H(M)$ is a strict one.
This explains the term ``weak $\Sl$-comodule'' for~$(C,M,u)$.

Any right $\Sl$-comodule~$M$ canonically gives
the weak $\Sl$-comodule~$(\Sl,M,u_P)$ where $u_P$ is the canonical
twisting cochain~\eqref{definition-canonical-twisting-cochain}.

A \newterm{morphism of weak $\Sl$-comodules} is a morphism of
right $\Ll$-modules~$f\colon M\otimes_u\Ll\to M'\otimes_{u'}\Ll$,
which we may write in the form
\begin{equation}
  f(m\otimes a)=\sum_{\pi\subset[r]}\doublesign\pi m\,
                  f_\pi(m)\otimes x_\pi\wedge a
\end{equation}
for some~$f_\pi\in\Hom_{-\degree\pi}(M,M')$.
If $f_\pi=0$ for~$\pi\ne\emptyset$, then
$f$ (and, by abuse of language, also $f_\emptyset$) is called \newterm{strict}.
Homotopies between weak $\Sl$-comodules are defined analogously.
This gives us the category~$\wComod\Sl$ of weak $\Sl$-comodules and their
morphisms. (This notation indicates that weak $\Sl$-comodules are in fact
right comodules over the reduced bar construction~$B\Ll$ of~$\Ll$
because a twisting cochain~$u\colon C\to\Ll$ with~$\epsilon_{\Ll}u=0$
is tantamount to a morphism of
coalgebras~$C\to B\Ll$, \cf~\cite[p.~77]{Smirnov:01}.)

If $f\colon M\to M'$ is a morphism of weak $\Sl$-comodules,
then this imposes certain conditions on its components~$f_\pi$,
namely
\begin{equation}\label{morphism-weak-comodules}
  \forall\,\pi\subset[r]\quad
  d(f_\pi)(m)=\sum_{\smashsubstack{(\mu,\nu)\vdash\pi\\\mu\ne\emptyset}}
    \sign{(\mu,\nu)}\Bigl(\sign\nu\,f_\nu(\gamma_\mu\cap m)
      -\doublesign\mu\nu\,\gamma_\mu'\cap f_\nu(m)\Bigr).
\end{equation}
In particular, $d(f_\emptyset) = 0$ and
$
   d(f_i)(m)
   = f_\emptyset(\gamma_i\cap m)-\gamma'_i\cap f_\emptyset(m)
$
for~$i\in[r]$.
Therefore, $f_\emptyset\colon M\to M'$ is a chain map of complexes and
induces an $\Sl$-equivariant map in homology.
We define $H(f)\colon H(M)\to H(M')$ to be this map~$H(f_\emptyset)$.
If $h\colon M\to M'$ is a homotopy between two morphisms $f$~and~$f'$, then
one verifies similarly that $h_\emptyset$ is a homotopy
between $f_\emptyset$~and~$f'_\emptyset$. Hence $H(f)=H(f')$ in this case.

\subsection{Weak $\Ll$-modules}\label{weak-modules}

Though we will not need it in our topological applications, we also introduce
the category~$\wMod\Ll$ of weak $\Ll$-modules for the sake of completeness.
A \newterm{weak \textmd{(left)} $\Ll$-module}~$N$ is a triple~$(A,N,v)$
where $N$ is a left module over an algebra~$A$
and $v\in\Hom_{-1}(\Sl,A)$ a twisting cochain satisfying~$v\iota_{\Sl}=0$.
(Hence~$N$ it is a left module over the reduced cobar construction of~$\Sl$.)
Similar to the previous case, we
identify~$v=\sum_{0\ne\alpha\in\N^r}c_\alpha\otimes\xi^\alpha\in A\otimes\Su$.
(Note that we always have $\sign{c_\alpha}=-1$.)
The differential on~$\Sl\otimes_v N$ is of the form
\begin{equation}
  d(s\otimes n)=s\otimes d n
    +\sum_{\alpha\ne0}\xi^\alpha\cap s\otimes c_\alpha\cdot n.
\end{equation}
(This is well-defined because $\xi^\alpha\cap s$ vanishes
for almost all~$\xi^\alpha$.)
The twisting cochain condition now reads
\begin{equation}
  \forall\,0\ne\alpha\in\N^r\qquad
  d c_\alpha=\sum_{\beta+\gamma=\alpha}c_\beta\cdot c_\gamma,
\end{equation}
in particular $d c_i = 0$, $d c_{ii} = c_i\cdot c_i$,
and~$d c_{i j} = c_i\cdot c_j+c_j\cdot c_i$ for~$i\ne j$.
Here we have written $c_i$~and~$c_{ii}$ for~$c_\alpha$
with a single non-vanishing component $\alpha_i=1$~and~$\alpha_i=2$,
respectively, and similarly for~$c_{i j}$. Hence the homology~$H(N)$ of
~$(A,N,u)$ carries a well-defined $\Ll$-module structure
defined by~$x_i\cdot[n]=[c_i\cdot n]$.

A \newterm{morphism of weak $\Ll$-modules}~$(A,N,v)\to(A',N',v')$
is a morphism of $\Sl$-comodules~$g\colon\Sl\otimes_v N\to\Sl\otimes_{v'}N'$,
which we can write as
\begin{equation}
  g(s\otimes n)=\sum_{\alpha\in\N^r}\xi^\alpha\cap s\otimes g_\alpha(n).
\end{equation}
Then the lowest component~$g_0$ is a chain map, and from~$
  d(g_i)(n)=g_0(c_i\cdot n)-c_i'\cdot g_0(n)$
we see that $g_0$ induces a $\Ll$-equivariant map in homology.
As before, the latter depends only on the (suitably defined)
homotopy class of~$g$.

\subsection{Koszul functors}\label{Koszul-algebraic}

The definition of the Koszul functors
\[
  \alg t\colon\wMod\Ll\to\wComod\Sl
  \and
  \alg h\colon\wComod\Sl\to\wMod\Ll
\]
is now almost obvious; we only have to switch
between left and right structures:
For any weak $\Ll$-module~$(A,N,v)$ we set
$\alg t N=\Sl\otimes_v N$ with the right $\Sl$-comodule structure
\begin{equation}
  T_{\Sl,\Sl\otimes_v N}(\Delta_{\Sl}\otimes1)\colon
    \Sl\otimes_v N\to\Sl\otimes\Sl\otimes_v N\to(\Sl\otimes_v N)\otimes\Sl,
\end{equation}
which is well-defined because $\Sl$ is cocommutative. Since this map is equal
to~$(1\otimes T_{\Sl N})(\Delta_{\Sl}\otimes1)$, the corresponding left
$\Su$-module structure is given by
\begin{equation}
    \sigma\cdot(s\otimes n) = \sigma\cap s\otimes n
    = (1\otimes\sigma)\Delta_{\Sl}(s)\otimes n.
\end{equation}
The left $\Ll$-module structure on~$\alg h M=M\otimes_u\Ll$,
$(A,M,u)$ a weak $\Sl$-comodule, is analogously defined by
\begin{equation}
    a\cdot(m\otimes b) = \sign a\doublesign a{m\otimes b}\,m\otimes b\wedge a
    = \sign a\doublesign a m\,m\otimes a\wedge b.
\end{equation}

Note that $\alg t N$~and~$\alg h M$ always possess strict structures.
In other words, the Koszul functors actually map
to $\rComod{\Sl}$~and~$\lMod\Ll$, respectively.
Moreover, both $\alg h\alg t R$~and~$\alg t\alg h R$ essentially give the
Koszul complex~$\K$.

\begin{theorem}\label{Koszul-duality-algebraic}
  The Koszul functors form an adjoint pair~$(\alg h,\alg t)$.
  Moreover, they induce quasi-inverse equivalences of categories in homology
  if restricted to the full subcategories of (co)modules
  bounded from below.
\end{theorem}

\begin{proof}
  We define natural transformations
  $\id\to\alg h\alg t$~and~$\alg t\alg h\to\id$ by the morphisms
  \begin{subequations}\label{Koszul-quasi-isomorphisms}
  \begin{gather}
    \Delta_{\Sl}\otimes1\otimes\iota_{\Ll}\colon
      \Sl\otimes_v N
      \to\Sl\otimes_P\alg h\alg t N
      =\Sl\otimes_P\bigl((\Sl\otimes_v N)\otimes_P\Ll\bigr),\\
    \epsilon_{\Sl}\otimes1\otimes\mu_{\Ll}\colon
      \big(\Sl\otimes_P(M\otimes_u\Ll)\bigr)\otimes_P\Ll
      =\alg t\alg h M\otimes_P\Ll
      \to M\otimes_u\Ll.
  \end{gather}
  \end{subequations}
  It is readily verified that the compositions
  \begin{align*}
    \alg t N\to\alg t(\alg h\alg t N) &= \alg t\alg h(\alg t N)\to\alg t N,\\
    \alg h M\to\alg h\alg t(\alg h M) &= \alg h(\alg t\alg h M)\to\alg h M
  \end{align*}
  are the respective identity morphisms. This proves the first claim.

  To see that they induce quasi-inverse equivalences in homology
  for (co)modules bounded from below, we first note that
  their lowest order components are of the form
  \begin{align*}
    \iota_{\Sl}\otimes1\otimes\iota_{\Ll}\colon
      N\to\alg h\alg t N=(\Sl\otimes_v N)\otimes_P\Ll, \\
    \epsilon_{\Sl}\otimes1\otimes\epsilon_{\Ll}\colon
      \Sl\otimes_P(M\otimes_u\Ll)=\alg t\alg h M\to M.
  \end{align*}
  Now filter the complexes
  \begin{align*}
    (\alg h\alg t N)_n &=
      \bigoplus_{p+q=n\vphantom X}\bigoplus_{p'+p''=p}
        \Sl_{p'}\oplus N_q\otimes\Ll_{p''},\\
    (\alg t\alg h M)_n &=
      \bigoplus_{p+q=n\vphantom X}\bigoplus_{q'+q''=q}
        \Sl_{q'}\oplus M_p\oplus\Ll_{q''}
  \end{align*}
  by $p$-degree, and consider $N$~and~$M$ as concentrated
  in $p$-degree~$0$ and $q$-degree~$0$, respectively.
  Then the above maps are filtration-preserving.
  The $E^1$-term for~$\alg h\alg t N$ is essentially the tensor product
  of the Koszul complex~$\K$ and~$H(N)$, so that we have an isomorphism
  between the $E^2$-terms, which are both equal to~$H(N)$.
  In the second case the tensor product~$\K\otimes M$ appears already on
  the $E^0$-level and therefore the isomorphism on~$E^1$.
  (Here we have used the assumptions $\epsilon_{\Ll}u=0$~and~$v\iota_{\Sl}=0$.)
\end{proof}

By replacing $\alg t N$~and~$\alg h M$ by arbitrary strict (co)modules in
formulas~\eqref{Koszul-quasi-isomorphisms}, one sees that the Koszul functors
are quasi-inverses between the subcategories of strict (co)modules
which are bounded from below.
One can show that the connecting morphisms (which now go the other way round)
are actually homotopy equivalences in the category of complexes, \ie,
after forgetting the (co)module structures \cite[Thm~1.6.3]{Franz:01}.
(See \cite[\S II.4]{HusemollerMooreStasheff:74} for a related result.)
Since this does not require the complexes to be bounded,
it implies in particular that the Koszul functors induce equivalences
in homology
between the categories of all $\Ll$-modules and all $\Sl$-comodules.
I conjecture that an analogous statement holds in the context
of weak (co)modules.
A careful discussion of Koszul duality between categories of
unbounded (strict) modules can be found in~\cite{Floystad:00a}.

\begin{proposition}\label{Koszul-preserve-quasi-isomorphism}
  The Koszul functors preserve quasi-isomorphisms between weak (co)modules
  bounded from below.
\end{proposition}

\begin{proof}
  This is again a spectral sequence argument,
  analogous to~\cite[Sec.~9]{GoreskyKottwitzMacPherson:98}
\end{proof}

This implies that the
Koszul functors $\alg h$~and~$\alg t$ descend to the derived categories
of (co)modules bounded from below,
which are obtained by localising at all quasi-isomorphisms.
Note, however, that the categories $\wMod\Ll$~and $\wComod\Sl$ are not abelian,
in contrast to those considered 
in \cite{GoreskyKottwitzMacPherson:98}~or~\cite{BernsteinGelfandGelfand:78}:
Any morphism of free $\Ll$-modules (with trivial differentials)
whose kernel is not free gives an example of a morphism of weak $\Sl$-comodules
without kernel.
Similarly, any morphism of finitely generated free $\Su$-modules
(again~$d=0$) with non-free kernel gives after dualising
a morphism of weak $\Ll$-modules without cokernel.

\section{Simplicial Koszul duality}\label{simplicity}

\subsection{Notation}\label{notation-2}

We will work in the category of simplicial sets,
see \cite{May:68} or~\cite{Lamotke:68} for expositions.
Recall that a simplicial set is an $\N$-graded set~$X=(X_n)$ together with
face maps~$\newsymbol{\partial_i}\colon X_n\to X_{n-1}$,
$0\le i\le n$ (for positive~$n$) and
degeneracy maps~$\newsymbol{s_i}\colon X_n\to X_{n+1}$,
$0\le i\le n$ satisfying certain commutation relations.
The basic example to keep in mind is of course the simplicial set of
singular simplices in a topological space.
Here $\partial_i\sigma$ is the composition of the
singular $n$-simplex~$\sigma$
with the inclusion~$\Delta_{n-1}\to\Delta_n$ of the $i$-th facet,
and $s_i\sigma$ the composition with the projection~$\Delta_{n+1}\to\Delta_n$
identifying the $i$-th vertex with its successor.
(From this one can deduce the commutation relations.)
Since simplicial sets are purely combinatorial objects,
one has a much greater flexibility in constructing them
than in the case of topological spaces.
This will be crucial when we define the connecting
natural transformation~$\tcomp$ in the next chapter.

We use the term \newterm{space} to refer to a simplicial set.
We write $C(X)$~and~$C^*(X)$ for the normalised chain resp.~cochain complex
of the space~$X$ with coefficients in~$R$,
and $f_*$~and~$f^*$ for the (co)chain map
induced by a map of spaces~$f\colon X\to Y$.
(Recall that $C(X)$ is obtained from the non-normalised chain complex
by dividing out the subcomplex of all degenerate simplices~$s_i\sigma$.
The projection is a homotopy equivalence, \cf~\cite[Sec.~VIII.6]{MacLane:75}.)
Note that for $X=\newsymbol\oneptspace$~a one-point space the complex~$C(X)$
is canonically isomorphic to~$R$.
The simplicial interval is denoted by~$\Deltaone$; its $n$-simplices
are the weakly increasing sequence~$(x_0,\ldots,x_n)$ of zeroes and ones.

To simplify notation, we introduce the abbreviations
\[
  \newsymbol{\partial_i^j}=
    \partial_i\circ\partial_{i+1}\circ\dotsb\circ\partial_j,
  \quad
  \partial_i^{i-1}=\id,
  \and
  \newsymbol{s_\mu}=s_{i_q}\circ\dotsb\circ s_{i_1},
  \quad
  s_\emptyset=\id
\]
for~$i\le j$ and any set~$\mu=\{i_1<i_2<\dotsb<i_q\}\subset\N$,
and also $\tpartial$ for the last face map, \ie,
$\tpartial=\partial_n$ in degree~$n$.

\subsection{The Eilenberg--Mac\,Lane maps}\label{products}

The Eilenberg--Zilber theorem states that the complexes
$C(X\times Y)$~and~$C(X)\otimes C(Y)$ are naturally homotopy-equivalent
for every pair~$X$,~$Y$ of spaces.
Particularly nice homotopy equivalences have been given by
Eilenberg and Mac\,Lane~\cite{EilenbergMacLane:54a}.
More precisely, they introduced certain maps
\begin{align*}
  \shuffle=\shuffle_{X Y} \colon C(X)\otimes C(Y) &\to C(X\times Y),\\
  \AW=\AW_{X Y} \colon C(X\times Y) &\to C(X)\otimes C(Y),\\
  H=H_{X Y}\colon C(X\times Y) &\to C(X\times Y)
\end{align*}
such that
\begin{subequations}\label{Eilenberg-Zilber-relations}
\begin{align}
  \AW\shuffle &= 1, \label{AW-shuffle-1}\\
  \shuffle\AW-1 &= d(H), \label{homotopy-shuffle-AW}\\
  \AW H &= 0, \\
  H\shuffle &= 0, \label{homotopy-shuffle-0}\\
  H H &= 0.
\end{align}
\end{subequations}
These maps are defined on the non-normalised complexes,
but descend to the normalised complexes, where they enjoy the properties
listed above.

\begin{subequations}\label{eilenberg-zilber-maps}
The \newterm{shuffle map}~$\newsymbol{\shuffle}$
carries the chain~$x\otimes y$, $x\in X_m$,~$y\in Y_n$, to
\begin{equation}\label{shuffle}
  \shuffle(x\otimes y)=\sum_{\smashsubstack{(\mu,\nu)\vdash(m,n)}}
    \sign{(\mu,\nu)}(s_\nu x,s_\mu y),
\end{equation}
where the sum extends over all $(m,n)$-shuffles (see Section~\ref{notation}).
The \newterm{Alexander--Whitney map}~$\newsymbol{\AW}$
is defined by
\begin{equation}\label{AW}
  \AW(x,y)
    =\sum_{i=0}^n\partial_{i+1}^n x\otimes\partial_0^{i-1}y
    =\sum_{i=0}^n(\tpartial)^{n-i}x\otimes(\partial_0)^i y
\end{equation}
for~$(x,y)\in(X\times Y)_n$.
The following non-recursive definition of the chain homotopy~$H$
is due to Rubio and Morace, \cf~\cite[Sec.~2]{GonzalezReal:99}:
\begin{equation}
  H(x,y) = \sum_{\smashsubstack{0\le i<j\le n\\(\mu,\nu)\vdash(j-i,n-j)}}
             (-1)^i\sign{(\mu,\nu)}
               \bigl(s_{\nu+i+1}s_i\partial_{j+1}^n x,
                 s_{\mu+i+1}\partial_{i+1}^{j-1}y\bigr),
\end{equation}
where $s_{\nu+k}=s_{\nu_{n-j}+k}\circ\cdots\circ s_{\nu_1+k}$.
\end{subequations}

The shuffle map and the Alexander--Whitney map are associative,
the former is also commutative in the sense that
\begin{diagram*}
  C(X)\otimes C(Y) & \rTo^{\shuffle_{X Y}} & C(X\times Y) \\
  \dTo<{T_{C(X),C(Y)}} & & \dTo>{{\tau_{X Y}}_*} \\
  C(Y)\otimes C(X) & \rTo^{\shuffle_{Y X}} & C(Y\times X)
\end{diagram*}
commutes, where $\tau_{X Y}(x,y)=(y,x)$, and $T$ is the transposition
of factors introduced in Section~\ref{notation}.
(See for example \cite[\S\,II.4]{Shih:62}~and~\cite[Ex.~12.26]{Dold:80}.)

The normalised chain complex of a space~$X$ is naturally a coalgebra
with comultiplication~$\AW\Delta_*\colon  C(X)\to C(X)\otimes C(X)$
and augmentation~$C(X)\to C(\oneptspace)=R$,
where $\Delta=\Delta_X\colon X\to X\times X$ is the diagonal.
In particular, we have cup and cap products as given
by equations \eqref{definition-cup-product}~and~\eqref{definition-cap-product}.
We shall also need the cohomological \newterm{cross product}~%
$\alpha\times\beta=(\alpha\otimes\beta)\AW\in C^*(X\times Y)$
of~$\alpha\in C^*(X)$~and~$\beta\in C^*(Y)$.

\begin{proposition}[Eilenberg--Moore~{\cite[\S 17.6]{EilenbergMoore:66}}]%
  \label{shuffle-morph-coalgebras}
  The shuffle map is a morphism of coalgebras.
\end{proposition}


\begin{proposition}[Shih]\label{X-Y-Z}
  For all spaces~$X$,~$Y$, and~$Z$, the following diagrams commute:
  \begin{diagram*}
    C(X\times Y)\otimes C(Z) & \rTo^{\shuffle_{X\times Y,Z}}
      & C(X\times Y\times Z)\\
    \dTo<{\AW_{X Y}\otimes1} & & \dTo>{\AW_{X,Y\times Z}} \\
    C(X)\otimes C(Y)\otimes C(Z) & \rTo_{1\otimes\shuffle_{Y Z}}
      & C(X)\otimes C(Y\times Z),
  \end{diagram*}
  \begin{diagram*}
    C(X)\otimes C(Y\times Z) & \rTo^{\shuffle_{X,Y\times Z}}
      & C(X\times Y\times Z) \\
    \dTo<{1\otimes\AW_{Y Z}} & & \dTo>{\AW_{X\times Y,Z}} \\
    C(X)\otimes C(Y)\otimes C(Z) & \rTo_{\shuffle_{X Y}\otimes1}
      & C(X\times Y)\otimes C(Z),
  \end{diagram*}
\end{proposition}

\begin{proof}
  Actually only the second diagram appears in~\cite[\S\,II.4]{Shih:62}.
  Like the first one it can be deduced from
  \theoremref {shuffle-morph-coalgebras}.
\end{proof}

\subsection{The Steenrod map}\label{section-Steenrod}

The cup product is not commutative, but only homotopy commutative.
A particularly nice homotopy is given by the \cuponeproduct.
I call the underlying map the \newterm{Steenrod map}~%
$\newsymbol{\ST}=\ST_{X Y}\colon C(X\times Y)\to C(X)\otimes C(Y)$.
It is not as fundamental as the previous maps, but defined as the composition
of the
``\index{Alexander--Whitney map!commuted}commuted Alexander--Whitney map''
\[
  \AWtilde_{\!X Y}=T_{C(Y),C(X)}\AW_{Y X}{\tau_{X Y}}_*\colon
    C(X\times Y)\to C(X)\otimes C(Y)
\]
and the Eilenberg--Mac\,Lane chain homotopy, namely~%
$\ST_{X Y}=\AWtilde_{\!X Y}H_{X Y}$.
It carries the non-degenerate
simplex~$(x,y)\in(X\times Y)_n$ to
\begin{equation}\label{ST-explicit}
  \ST(x,y)=\sum_{\smashsubstack{0\le i<j\le n}}(-1)^{i+j+i j}\,
    \partial_0^{i-1}\partial_{j+1}^n x\otimes\partial_{i+1}^{j-1}y
\end{equation}
of degree~$n+1$, \cf~\cite[Thm.~3.1]{GonzalezReal:99}.

The \newterm[cross-one product@\protect\crossoneproduct]{\crossoneproduct}~%
$\newsymbol{\crossone}\colon C^*(X)\otimes C^*(Y)\to C^*(X\times Y)$
and the \newterm[cup-one product@\protect\cuponeproduct]{\cuponeproduct}~%
$\newsymbol{\cupone}\colon C^*(X)\otimes C^*(X)\to C^*(X)$
are (up to sign) derived from the Steenrod map
like cross and cup product from the Alexander--Whitney map, \ie,
\begin{align*}
  \alpha\crossone\beta &= \sign\alpha\sign\beta(\alpha\otimes\beta)\ST\colon
    C(X\times Y)\to R\otimes R=R, \\
  \alpha\cupone\beta &= (\alpha\crossone\beta)\Delta_*.
\end{align*}

\begin{proposition} $ $
  \begin{enumerate}
  \item\label{cross1-zero}
    A {\crossoneproduct} or {\cuponeproduct} vanishes if one factor
    is of degree zero.
  \item\label{cross1-shuffle}
    The {\crossoneproduct} vanishes on the image of the shuffle map for all
    pairs of spaces.
  \item The {\cuponeproduct} is a homotopy between the cup product
    and the cup product with commuted factors:
    \[
      d(\alpha\cupone\beta) =
        \alpha\cup\beta
        -\doublesign\alpha\beta\,\beta\cup\alpha
        -d\alpha\cupone\beta
        -\sign a\,\alpha\cupone d\beta.
    \]
  \item It is also a \emph{left} derivation of the cup product
    (Hirsch formula):
    \[
      \alpha\cupone(\beta\cup\gamma) =
        (\alpha\cupone\beta)\cup\gamma
        +\doublesign\alpha\beta\sign\beta\,\beta\cup(\alpha\cupone\gamma).
    \]
  \end{enumerate}
  Here $\alpha$,~$\beta$, and~$\gamma$ denote cochains on some space.
\end{proposition}

One usually defines the {\cuponeproduct} such that it becomes
a right derivation, but this would be less convenient in our applications.
We remark in passing that the signs in the above formulas are as predicted by
the sign rule if we write the {\cuponeproduct}
as map~$\cupone(\alpha\otimes\beta)$.

\begin{proof}
  The first assertion follows directly from the explicit
  formula~\eqref{ST-explicit} and
  the second and third
  from equations~\eqref{Eilenberg-Zilber-relations}.
  The Hirsch formula requires a lengthy calculation.
\end{proof}

\begin{proposition}\label{properties-Steenrod}
  Let $X$, $Y$, and~$Z$ be spaces.
  \begin{enumerate}
  \item The following diagram commutes:\label{Steenrod-shuffle-associative}
    \begin{diagram*}
      C(X)\otimes C(Y\times Z) & \rTo^{\shuffle_{X,Y\times Z}}
        & C(X\times Y\times Z) \\
      \dTo<{1\otimes\ST_{Y Z}} & & \dTo>{\ST_{X\times Y,Z}} \\
      C(X)\otimes C(Y)\otimes C(Z) & \rTo_{\shuffle_{X Y}\otimes1}
        & C(X\times Y)\otimes C(Z).
    \end{diagram*}
  \item For all~$w\in C(X\times Y)$ and~$z\in C(Z)$ with~$\degree z\le1$
    one has\label{Steenrod-shuffle-associative-1}
    \[
      \ST_{X,Y\times Z}\shuffle_{X\times Y,Z}(w\otimes z)
      = (1\otimes\shuffle_{Y Z})(\ST_{X Y}\otimes1)(w\otimes z).
    \]
  \end{enumerate}
\end{proposition}

These are
are analogues of \theoremref{X-Y-Z} for the Steenrod map,
but this correspondence is only partial because
part~(\ref{Steenrod-shuffle-associative-1}) above is \emph{false}
for general~$z$.

\begin{proof}
  This is verified by direct calculations, see~\cite[Appendix~8]{Franz:01}.
  The first part can alternatively be deduced from a large commutative
  diagram involving~\cite[\S\,II.4, Lemme~3]{Shih:62}.
\end{proof}

\subsection{Groups and group actions}\label{groups}

If not specified otherwise, a \newterm{group}~$G$ will mean a simplicial group
\ie, a simplicial object in the category of groups (or, equivalently,
a group object in the category of simplicial sets), \cf~\cite[\S 17]{May:68}.
When we are careful about notation, we write $\newsymbol{1_n}$ for the
unit element of the set-theoretic group~$G_n$ of $n$-simplices of~$G$.
For any group~$G$, we denote the category of left $G$-spaces
by~$\newsymbol{\Space G}$.

Let $\mu\colon G\times G\to G$ be the multiplication of~$G$,
$\lambda\colon G\to G$~the inversion and
$\iota\colon1\to G$~the unit.
The chain complex~$C(G)$ is an algebra with
Pontryagin product~$\mu_*\shuffle_{GG}$ and unit~$\iota_*$.
It follows from \theoremref{shuffle-morph-coalgebras} that $C(G)$ is actually
an (associative and coassociative) Hopf algebra.
By commutativity of the shuffle map, $C(G)$ is commutative if $G$ is.

If $X$ is a left (or right) $G$-space,
then $C(X)$ is canonically a left (or right) $C(G)$-module.
(Use again the shuffle map.)
If the action of~$G$ on~$X$ is trivial, then so is that of~$C(G)$ on~$C(X)$
because it factors through the augmentation~$C(G)\to C(1)=R$.
One can switch between left and right actions
by defining~$x\cdot g=g^{-1}\cdot x$ for $x\in X$~and~$g\in G$.
The two corresponding $C(G)$-module structures on~$C(X)$ are related by
\begin{equation}\label{comparison-C-G-module-structures}
  c\cdot a=\doublesign c a\,\lambda_*(a)\cdot c
\end{equation}
for $a\in C(G)$~and~$c\in C(X)$.
This again follows from the shuffle map's commutativity.

\begin{proposition}
  Let $G$ be a group and $X$~and~$Y$ both left or both right $G$-spaces.
  \begin{enumerate}
  \item The\label{shuffle-equivariant} shuffle map~$\shuffle_{X Y}$ is
    $C(G)$-equivariant if $G$ operates trivially on either space.
  \item The\label{AW-equivariant}
    Alexander--Whitney map~$\AW_{X Y} 
    $ is $C(G)$-equivariant.
  \item The\label{ST-equivariant}
    Steenrod map~$\ST_{X Y} 
    $ is $C(G)$-equivariant if $G$ operates trivially on~$Y$.
  \item It\label{ST-equivariant-1}
    is also equivariant with respect to an~$a\in C(G)$ if
    $\degree a\le1$ and $G$ operates trivially on~$X$.
  \item The\label{C-X-C-G-coalgebra}
    chain complex~$C(X)$ is a $C(G)$-coalgebra.
  \end{enumerate}
\end{proposition}

(See Section~\ref{notation} for the definition of a $C(G)$-coalgebra.)

\begin{proof}
  Notice first that it is not important whether $G$ acts from the left or
  from the right because we may always pass from one to the other
  by redefining the action and then,
  using equation~\eqref{comparison-C-G-module-structures},
  go back to the original one on the chain level.

  The first assertion is a consequence of the properties of the shuffle map.
  The non-trivial part of the last claim, the equivariance of
  the comultiplication, follows from the second assertion.
  This in turn is once again a consequence of
  \theoremref{shuffle-morph-coalgebras}.
  In the same way \theoremref{properties-Steenrod}
  proves the third and fourth claim.
\end{proof}

\begin{proposition}\label{chain-functor-Space-2-Mod}
  The chain functor~$C$ is a well-defined
  homotopy-preserving functor~$\Space G\to\lMod{C(G)}$.
\end{proposition}

\begin{proof}
  This is a routine verification using the properties
  of the Eilenberg--Mac\,Lane maps
  as given in Section~\ref{products}.
\end{proof}

\subsection{Spaces over a base space}

Let $B$ be a space. A \newterm[space!over a base space]{space over~$B$} is
a map of spaces~$p\colon Y\to B$.
We will usually refer to~$p$ by~$Y$, the map~$p$ being understood.
A \newterm{morphism~$f\colon p\to p'$ of spaces
over~$B$} is a map~$f\colon Y\to Y'$ such that~$p'\circ f=p$.
Similarly, homotopies between spaces over~$B$
are homotopies~$Y\times\Deltaone\to Y'$ commuting
with the projections to~$B$.
(Here $\Deltaone$ is the simplicial interval, \cf~Section~\ref{notation-2}.)
We denote the category of spaces over~$B$ by~$\newsymbol{\Spaceover B}$.

If $Y\to B$ is a space over~$B$, then $C(Y)$ is a right $C(B)$-comodule
via the map~$\AW_{Y B}
\left.\tilde\Delta_Y\right._{\!*}
\colon C(Y)\to C(Y)\otimes C(B)$, where
$\Deltap_Y$ denotes the canonical map of spaces~$Y\to Y\times B$.

\begin{proposition}\label{chain-functor-Spaceover-2-Comod}
  The chain functor~$C$ is a well-defined
  homotopy-preserving functor~$\Spaceover B\to\rComod{C(B)}$.
\end{proposition}

\begin{proof}
  This is again routine. The verification that one gets
  equivariant chain homotopies uses \theoremref{X-Y-Z}.
\end{proof}

\subsection{Fibre bundles}\label{fibre-bundles}

Our definition of a simplicial fibre bundle~\cite{BarrattGugenheimMoore:59}
is not exactly that of~\cite[\S 18]{May:68}:
Let $B$ be a space and $G$ a group. A \newterm{twisting function} for~$B$
with values in~$G$ is a map of graded sets~$\tau\colon B_{>0}\to G$
of degree~$-1$ such that for~$b\in B_n$ and~$i<n$
\begin{xalignat*}2
  \tau(\partial_i b) &= \partial_i\tau(b), &
  \tau(\partial_n b)
    &= \tau(\partial_{n-1}b)\bigl(\partial_{n-1}\tau(b)\bigr)^{-1},\\
  \tau(s_i b) &= s_i\tau(b), &
  \tau(s_n b) &= 1_n
\end{xalignat*}
Given a left $G$-space~$F$, one defines the
twisted Cartesian product~$B\times_\tau F$,
whose only difference from~$B\times F$ lies in the last face map, which is now
\[
  \tpartial(b,f) = \bigl(\tpartial b,\tau(b)\tpartial f\bigr).
\]
The corresponding \newterm{fibre bundle} is the obvious
projection~$B\times_\tau F\to B$.

We shall need the simplicial version of the
Leray--Serre spectral sequence of a fibre bundle:
Filter $C(B\times_\tau F)$ by the $p$-skeletons of the base~$B$.
More precisely, a non-degenerate simplex~$(f,b)\in C_n(B\times_\tau F)$
belongs to~$F_p C(B\times_\tau F)$ if $b$ is (at least)
$(n-p)$-fold degenerate. This leads to a spectral sequence,
natural in $B$~and~$F$, whose first terms are for connected~$G$ equal to
\begin{subequations}
\begin{align}
  E^1_{p q}(B, F) &= C_p(B;H_q(F)),\\
  E^2_{p q}(B, F) &= H_p(B;H_q(F))
\end{align}
\end{subequations}
as comodules over $C(B)$~and~$H(B)$, respectively,
and also as $H(G)$-modules in the case~$F=G$.\footnote{More precisely,
these terms arise from a filtration of a certain
twisted tensor product~$C(B)\otimes_u C(F)$, which is related
to~$C(B\times_\tau F)$ by a filtration-preserving homotopy-equivalence.
See the proof of \theoremref{Leray-Serre-allowable-base} for more details
and references.}
An analogous spectral sequence exist in cohomology with
\begin{equation}\label{Leray-Serre-cohomological}
  E_2^{p q}(B, F) = H^p(B;H^q(F)).
\end{equation}

\subsection{Classifying spaces and universal bundles}\label{classifying-spaces}

In order to define Koszul functors in the simplicial setting, we need the
simplicial construction of universal bundles and classifying spaces.
Our definitions again differ from~\cite[\S 21]{May:68}:

For any group~$G$, the \newterm{classifying space}~$\newsymbol{BG}$
is the space with set of $n$-simplices
\[
  BG_n=G_0\times\dotsb\times G_{n-1}
\]
for~$n\in\N$. We write the simplices of~$BG$ in the
form
\[
  [g_0,\dotsc,g_{n-1}]\in BG_n,
  \qquad\text{also}\qquad
  \newsymbol{b_0}:=[]\in BG_0
\]
for the unique vertex of~$BG$.
The face and degeneracy maps are given by
\begin{align*}
  \partial_i[g_0,\dotsc,g_{n-1}] &=
    [g_0,\dotsc,g_{i-2},g_{i-1}\partial_0g_i,
    \partial_1g_{i+1},\dotsc,\partial_{n-i-1}g_{n-1}],\\
  s_i[g_0,\dotsc,g_{n-1}] &=
    [g_0,\dotsc,g_{i-2},g_{i-1},1_i,s_0g_i,
    s_1g_{i+1},\dotsc,s_{n-i-1}g_{n-1}].
\end{align*}
(Undefined components, such as~$\partial_{-1}g_{n-1}$, are supposed to be
omitted when applying these formulas for given values of~$i$ and~$n$.)
The map of graded sets
\[
  \tau_{BG}\colon BG_{>0}\to G,
  \qquad
  [g_0,\dotsc,g_{n-1}]\mapsto g_{n-1}
\]
is a twisting function for~$BG$.
The principal bundle~$EG=BG\times_{\tau_{BG}}G\to BG$ is called
the universal $G$-bundle.
$G$ acts freely from the right on its total space.
The construction of classifying spaces and universal bundles
is functorial and compatible with products, \ie,
\begin{equation}\label{BG-EG-products}
  B(G\times H)=BG\times BH
  \and
  E(G\times H)=EG\times EH
\end{equation}
for any pair~$G$,~$H$ of groups.

\medbreak

Another way to look at~$EG$ is the following: We have $EG_n=BG_{n+1}$ as sets,
moreover $\partial^{EG}_i=\partial^{BG}_i$ for~$i\le n$
and likewise for degeneracy maps. This implies that $EG$ is homotopy-equivalent
to~$s_0BG_0\subset BG_1=EG_0$, \ie, that it is (non-equivariantly) contractible
to~$e_0$.
(This will be proven along the way in \theoremref{simplicial-composition}.)
The last degeneracy map on~$BG$ now gives rise to the following map~$\pS$
of degree~$1$, \cf~\cite[Def.~21.1]{May:68}:
\begin{equation}
  \newsymbol{\pS}=\pS_G\colon EG\to EG,
  \qquad
  ([g_0,\dotsc,g_{n-1}],g_n)
  \mapsto
  ([g_0,\dotsc,g_{n-1},g_n],1_{n+1}).
\end{equation}
Rewriting the commutation relations in terms of~$\pS$, we see that it satisfies
for all~$e\in EG_n$ and~$0\le i\le n$ the identities
\begin{subequations}\label{d-S-and-s-S}
  \begin{align}
    \partial_i \pS e &= \begin{cases}
       \pS\partial_i e & \text{if~$n>0$},\\
       e_0 & \text{if~$n=0$},
    \end{cases}
    &
    \partial_{n+1}\pS e &= e, \\
    s_i \pS e &= \pS s_i e, & s_{n+1} \pS e &= \pS \pS e.
  \end{align}
\end{subequations}
In particular, $\pS$ passes to a map~$C(EG)\to C(EG)$
on the normalised complex.
It will be convenient to rescale it and to define
\begin{equation}
  \newsymbol{S}\colon C(ET)\to C(ET),
  \quad
  e\mapsto-\sign e\,\pS e.
\end{equation}
The composition~$S\circ S$ vanishes, and~$S e_0=-s_0 e_0=0$.
Moreover, $S$ is a homotopy between the identity and the chain map induced
by the retraction~$EG\to e_0$:
\begin{equation}\label{S-homotopy}
  S d e+d S e=\begin{cases}
    e & \text{if~$\degree e>0$},\\
    e-e_0 & \text{if~$\degree e=0$}
  \end{cases}
\end{equation}
for all non-degenerate~$e\in EG$.
(The map~$S$ is actually the chain homotopy induced
by the aforementioned strong deformation retraction of~$EG$ to~$e_0$,
see the proof of \theoremref{simplicial-composition}).

Note that $S$ is compatible with products, \ie, $S_{G\times H}=S_G\times S_H$.
We need to know how $S$ interacts
with the Eilenberg--Mac\,Lane maps.

\begin{proposition}\label{cone-homotopy}
  The following identities hold for all groups $G$~and~$H$:
  \begin{subequations}
  \begin{align}
    \label{shuffle-S}
    \shuffle_{EG,EH}(S_G\otimes S_H)
      &= S_{G\times H}\shuffle_{EG,EH}(1\otimes S_H-S_G\otimes1),\\
    \AW_{EG,EH}S_{G\times H}(e,e')
      &= S_G(e)\otimes e_0'+(1\otimes S_H)\AW_{EG,EH}(e,e') \\
  \intertext{for~$(e,e')\in EG\times EH$ and~$e_0'$ the canonical
    basepoint of~$EH$,}
    \label{ST-S}
    \ST_{EG,EH}S_{G\times H}
      &= (S_G\otimes S_H)\AW_{EG,EH}-(1\otimes S_H)\ST_{EG,EH}.
  \end{align}
  \end{subequations}
\end{proposition}

\begin{proof}
  These are direct computations, see~\cite[Appendix~9]{Franz:01}.
  (There $BG$,~$EG$, and $S$ have slightly different definitions.)
  For~\eqref{shuffle-S} one splits up the sum over all
  $(m+1,n+1)$-shuffles~$(\mu,\nu)$ depending on whether the maximum~$m+n+1$
  is contained in~$\mu$~or in~$\nu$.
\end{proof}

\subsection{Simplicial Koszul functors}\label{simplicial-functors}

We are now in the position to introduce, for any group~$G$, the
simplicial Koszul functors
between the categories of $G$-spaces and spaces over~$BG$.

Let $X$ be a left $G$-space. The
\newterm{Borel construction}
\begin{equation}
   \newsymbol{EG\timesunder G X}=BG\times_{\tau_{BG}}X\to BG 
\end{equation}
is a space over~$BG$.
The map of spaces
\begin{equation}
  \newsymbol{q_X}\colon EG\times X\to EG\timesunder G X,
  \qquad
  \bigl((b,g),x\bigr)\mapsto(b,g x),
\end{equation}
is the quotient of~$EG\times X$
by the $G$-action~$g\cdot(e,x)=(e g^{-1},g x)$.
We record the following observation, which will be used in
Section~\ref{comparing-functors}:

\begin{lemma}\label{q-X-shuffle-equivariant}
  The composition
  \[
    {q_X}_*\shuffle_{EG,X}\colon C(EG)\otimes C(X)
      \to C(EG\times X)\to C(EG\timesunder G X)
  \]
  is a morphism of right $C(BG)$-comodules.
\end{lemma}

Here $C(EG\timesunder G X)$ is a right comodule because $EG\timesunder G X$
is a space over~$BG$, and $C(EG)\otimes C(X)$ is one because $C(EG)$ is,
\ie, the structure map is
\[
  (1\otimes T_{C(BG),C(X)})(\Delta_{C(EG)}\otimes1)\colon
    C(EG)\otimes C(X)\to \bigl(C(EG)\otimes C(X)\bigr)\otimes C(BG).
\]

\begin{proof}
  The map~${q_X}_*\shuffle$ can be written as composition
  \[
    C(EG)\otimes C(X)\to C(X)\otimes C(EG)\to C(X\times EG)
      \to C(X\timesunder G EG)\to C(EG\timesunder G X)
  \]
  by the shuffle map's commutativity. Hence it suffices to prove
  the claim for the map~$C(X)\otimes C(EG)\to C(X\times EG)$.
  Here equivariance follows from the commutative diagram
  \begin{diagram*}
    C(X)\otimes C(EG) & \rTo^{1\otimes \Deltap_*}
      & C(X)\otimes C(EG\times BG)
      & \rTo^{1\otimes \AW} & C(X)\otimes C(EG)\otimes C(BG) \\
    \dTo>\shuffle & & \dTo>\shuffle & & \dTo>{\shuffle\otimes1} \\
    C(X\times EG) & \rTo^{(\id,\Deltap)_*}
      & C(X\times EG\times BG)
      & \rTo^{\AW} & C(X\times EG)\otimes C(BG) \\
  \end{diagram*}
\end{proof}

The simplicial Koszul functor
\[
  \newsymbol{\Sim t}\colon\Space G\to\Spaceover{BG}
\]
assigns to each left $G$-space~$X$ the space~$EG\timesunder G X$
(more precisely, the projection~$EG\timesunder G X\to BG$\,),
and to each morphism~$X\to X'$
the induced morphism~$EG\timesunder G X\to EG\timesunder G X'$.

Given a map~$p\colon Y\to BG$, we can pull back the universal $G$-bundle
to obtain the right $G$-space~$Y\times_{\tau_{BG}\circ p}G$,
which we also write as
$\newsymbol{Y\timesover{BG}EG}$ 
since it is the fibre product of $Y$~and~$EG$ over~$BG$.
The Koszul functor
\[
  \newsymbol{\Sim h}\colon\Spaceover{BG}\to\Space G
\]
assigns to each space~$Y$ over~$BG$ the space~$Y\timesover{BG}EG$
with the opposite, \ie, \emph{left} $G$-action,
and to each morphism~$Y\to Y'$ the induced
morphism~$Y\timesover{BG}EG\to Y'\timesover{BG}EG$.

Both functors preserve homotopies.

\medbreak

Let $X$ be a left $T$-space and $p\colon Y\to BG$ a space over~$BG$.
In order to introduce certain natural transformations
\begin{equation}
  \Sim P\colon\Sim h\Sim t\to\id
  \and
  \Sim I\colon\id\to\Sim t\Sim h,
\end{equation}
we use the isomorphisms of graded sets $\Sim h\Sim t X\cong BG\times X\times G$
and $\Sim t\Sim h Y\cong BG\times Y\times G$ and define the morphisms
\begin{subequations}\label{definitions-P-I}
\begin{gather}
  \Sim P_X\colon\Sim h\Sim t X\to X,
  \qquad
  (b,x,g)\mapsto g^{-1}x,\\
  \Sim I_Y\colon Y\mapsto\Sim t\Sim h Y,
  \qquad
  y\mapsto(p(y),y,1)
\end{gather}
\end{subequations}
for $X\in\Space G$~and~$Y\in\Spaceover{BG}$, respectively.
(That the above formulas define morphisms in these categories
will be become evident in the proof of \theoremref{simplicial-composition}.)

\begin{proposition}
  The simplicial Koszul functors form an adjoint pair~$(\Sim h,\Sim t)$.
\end{proposition}

\begin{proof}
  A look at the definitions of $\Sim I_Y$~and~$\Sim P_X$
  shows that the compositions
  \begin{gather*}
    \Sim h Y\longrightarrow\Sim h(\Sim t\Sim h Y)
      =\Sim h\Sim t(\Sim h Y)\longrightarrow\Sim h Y, \\
    \Sim t X\longrightarrow\Sim t\Sim h(\Sim t X)
      =\Sim t(\Sim h\Sim t X)\longrightarrow\Sim t X.
  \end{gather*}
  are the respective identities.
\end{proof}

\begin{theorem}\label{simplicial-composition}
  The morphisms $\Sim P_X$~and~$\Sim I_Y$ are homotopy equivalences of spaces
  for all $G$-spaces~$X$ and all spaces~$Y$ over~$BG$.
\end{theorem}

This is a (partial) analogue of~\cite[Remark~1.7]{AlldayPuppe:99}
in the simplicial setting.
It parallels the duality between the algebraic
Koszul functors, \cf~in particular the comments made
after \theoremref{Koszul-duality-algebraic}.
Analogously to that case,
it implies that the simplicial Koszul functors become quasi-inverse
equivalences of categories after localising the categories
$\Space G$~and~$\Spaceover{BG}$ with respect to all morphisms
which are homotopy equivalences of spaces.
The present analysis will prove useful in Section~\ref{main-allowable}.

\begin{proof}
   We start with~$\Sim I_Y$ and claim that the map of spaces
   \[
     \Sim J_Y\colon\Sim t\Sim h Y\to Y,
     \qquad
     (b,y,g)\mapsto y,
   \]
   is a homotopy inverse.
   We clearly have~$\Sim J_Y\Sim I_Y=\id_Y$.
   In the representation used above, the last face map
   of~$\Sim t\Sim h Y$ is of the form
  \[
    \tpartial(b,y,g)=\bigl(
      \tpartial b,\tpartial y,\tau_Y(y)(\tpartial g)\tau_{BG}(b)^{-1}
    \bigr)
  \]
  with~$\tau_Y=\tau_{BG}\circ p$.
  From this it follows that $\Sim J_Y$ is a map of spaces and $\Sim I_Y$ a map
  over~$BG$. (Recall that $\Sim t\Sim h Y$ is a space over~$BG$
  via projection onto the $b$-coordinate in our representation.)
  It is convenient to reorder the factors
  and to apply the group inversion to~$G$ so that
  \[
    \tpartial(y,b,g)=\bigl(
      \tpartial y,\tpartial b,\tau_{BG}(b)(\tpartial g)\tau_Y(y)^{-1}
    \bigr)
  \]
  This shows that $\Sim J_Y$ is essentially a bundle with base~$Y$,
  twisting function~$\tau_Y$, and fibre~$EG$
  (with left $G$-action~$g\cdot e=e\cdot g^{-1}$).
  We may therefore write the last face map as
  \[
    \tpartial(y,e)=\bigl(\tpartial y,\tpartial e\cdot\tau_Y(y)^{-1}\bigr).
  \]
  Note that for~$Y=\oneptspace$ the total space is just~$EG$.
  Hence we will now in particular show that $EG$ is contractible,
  as announced in the preceding section.

  In order to define
  a homotopy~$h\colon\Sim t\Sim h Y\times\Deltaone\to\Sim t\Sim h Y$
  from~$\Sim I_Y\circ \Sim J_Y$ to the identity of~$Y$,
  we need the map~$\pS$ introduced in Section~\ref{classifying-spaces},
  which was the last degeneracy map of~$BG$ carried over to~$EG$.
  Now $h$ is recursively given by
  \[
    h(y,e,x)=\begin{cases}
      (y,(s_0)^{\degree y}e_0) & \text{if~$x=(0)$}, \\
      \bigl(y,\pS(e'\cdot\tau_Y(y))\bigr) &
        \text{if~$x_0=0$, but~$x\ne(0)$}, \\
      (y,e) & \text{if~$x_0=1$},
    \end{cases}
  \]
  where $e'$ is determined by~$(\tpartial y,e')=h(\tpartial(y,e,x))$,
  and $x_0$ denotes the leading element of the sequence~$x\in\Deltaone$.
  (Recall that the simplices in~$\Deltaone$ are the weakly increasing
  sequences composed of zeroes and ones.)
  The verification that $h$ is a homotopy as claimed is elementary,
  but somewhat lengthy.
  (It can essentially be found in~\cite[Appendix~10]{Franz:01}.)

  For $Y=\oneptspace$, an~$e\in EG_n$ and $x$ consisting of $k$~zeroes followed
  by $n+1-k$~ones the formula for~$h$
  simplifies to~$h(e,x)=\pS^k\tpartial^k e$
  (which is to be read as~$e_0$ for~$k=n+1$).
  The induced chain homotopy on~$C(EG)$ is just~$S$.

  \goodbreak
 
  We now consider the transformation~$\Sim P$. Let $X$ be a $T$-space.
  That $\Sim P_X$ and
  \[
    \Sim Q_X\colon X\to\Sim h\Sim t X,
    \qquad
    x\mapsto\bigl((s_0)^{\degree x}b_0,x,1)
  \]
  are maps of spaces follows from the formula
  \begin{equation*}
    \tpartial(b,x,g)
      =\bigl(\tpartial b,\tau_{BG}(b)\tpartial x,\tau_{BG}(b)\tpartial g\bigr)
  \end{equation*}
  for the last face map of~$\Sim h\Sim t X$.
  The map
  \begin{equation*}
    \Sim h\Sim t X \to EG\times X,
    \qquad
    (b,x,g)\mapsto (b,g,g^{-1}x)
  \end{equation*}
  is an isomorphism of left $G$-spaces.
  Here $G$ acts on~$EG\times X$ by~$g(e,x)=(e g^{-1},g x)$.
  Under this isomorphism the map~$\Sim P_X$ corresponds to the
  canonical $G$-equivariant projection onto~$X$
  and $\Sim Q_X$ to the inclusion of~$X$ over~$e_0$ in~$EG\times X$.
  Since $EG$ is contractible, as just shown, these maps
  are homotopy-equivalences.
\end{proof}

\section{Comparing the functors}\label{main-chapter}

\subsection{Tori}\label{tori}

From now on we focus on tori, which are denoted by the letter~$T$
instead of the~$G$ used so far for groups.
We start with our definition of circles and tori.
Readers who have left out Section~\ref{classifying-spaces} can think
of compact tori~$\cong(S^1)^r$ or algebraic tori~$\cong(\C^*)^r$
and should skip the next paragraph.

A \newterm{simplicial circle} is a group isomorphic to the
classifying space~$B\Z$, which is a group by componentwise multiplication
because $\Z$ is commutative.
More prosaically, it is isomorphic to the subgroup of the group of
singular simplices in~$S^1$ generated by a loop at~$1$.
We define a \newterm{circle} as a group containing a simplicial circle as a
subgroup such that the inclusion is a quasi-isomorphism.
We let~$S^1$ denote any circle, and a ``loop at~$1\in S^1$''
refers to a generator of the embedded simplicial circle.
A \newterm{torus} of
\newterm[rank!of a torus]{rank}~$r$ is a group isomorphic to an $r$-fold
product of circles.

Given a torus~$T$, we fix once and for all a decomposition
into circles. We write $x'_i\in C(T)$ for a loop at~$1\in T$
around the $i$-th factor~$S^1$ (or~$\C^*$),
and similarly for its homology class~$x_i$.
Moreover, we now take $H_{\minusdegd}(T)$ as the free
$R$-module~$P$ (concentrated in degree~$\minusdegd$)
that was the starting point for the algebraic constructions
in Section~\ref{Koszul-complex}.
Note that we have a canonical isomorphism of Hopf algebras~$H(T)=\Ll$
as well as one of coalgebras~$H(BT)=\Sl$.

The following lemma\notheoremref[lemma]{section-Hopf-algebras} is the reason
why we concentrate on tori in this paper.\footnote{As mentioned
in~\cite[Sec.~12]{GoreskyKottwitzMacPherson:98}, one can construct
a morphism of algebras~$H(G)\to C(G)$ as in the
lemma\notheoremref[lemma]{section-Hopf-algebras}
whenever $H(G)=\Ll$ is an exterior algebra and the generators~$x_i\in\Ll$
can be represented by conjugation-invariant cycles~$x_i'\in C(G)$.
But contrary to the claim in~\cite{GoreskyKottwitzMacPherson:98},
this is not possible in general because all singular simplices appearing in a
conjugation-invariant chain necessarily map to the centre of~$G$.
The example $G=SU(3)$ shows that the use of subanalytic chains is no remedy:
Apart from the finite centre, all conjugation classes have dimension $4$
or~$6$. Hence there can be no conjugation-invariant subanalytic set
supporting a representative of the degree~$3$ generator.}

\begin{lemma}\label{section-Hopf-algebras}
  The assignment~$\Ll\to C(T)$,
  $x_{i_1}\wedge\dotsb\wedge x_{i_q}\mapsto x'_{i_1}\dotsm x'_{i_q}$
  is a quasi-isomorphism of Hopf algebras.
\end{lemma}

\begin{proof}
  The map is well-defined because $x_i'\cdot x_j'=-x_j'\cdot x_i'$
  by the shuffle map's commutativity. Since we have chosen each $x_i'$ to be
  a loop at~$1$, its image under~$\Delta_{C(T)}$ is
  $x_i'\otimes1+1\otimes x_i'$.
  Together with \theoremref{shuffle-morph-coalgebras} this shows
  that the coalgebra structures
  are compatible. It is clear that the map is a quasi-isomorphism.
\end{proof}

This map gives us a homotopy-preserving functor
from the category of modules over~$C(T)$
to those over~$\Ll$. In order to avoid a too clumsy notation,
we incorporate it into the chain functor.
Hence $\newsymbol{C}(X)\in\lMod\Ll$ for~$X\in\Space T$.
In~\cite{GoreskyKottwitzMacPherson:98}
the resulting $\Ll$-module structure on chain and cochain
complexes of $T$-spaces is called the ``sweep action''. 

\medbreak

To analogously view $C(BT)$-comodules as $\Sl$-comodules, we would need a
quasi-isomorphism of coalgebras~$C(BT)\to\Sl$. This does not exist
for~$r\ge2$, essentially because $C(BT)$ is not cocommutative.
But $C(BT)$-comodules can be turned into weak $\Sl$-comodules:

\begin{proposition}[Gugenheim--May]\label{functor-weak-S}
  Any choice of representatives~$\xi_i'\in C^*(BT)$
  of a set of generators~$\xi_i\in\Sl^2=H^2(BT)$
  canonically defines a twisting cochain
  \[
    t=\sum_{\emptyset\ne\pi\subset[r]}x_\pi\otimes\xi'_\pi
    \in\Hom_{-1}(C(BT),\Ll)
    \qquad\hbox{with}\qquad
    \epsilon_{\Ll}t=0
  \]
  by recursively setting
  \[
    \xi'_\pi=\xi'_{\pi^+}\cupone\xi'_{\pi'},
  \]
  where $\pi^+$ is the maximum of~$\pi$ and $\pi'\ne\emptyset$ the other
  elements.
\end{proposition}

Since this result is purely algebraic, it holds for any group~$G$ and any
coefficient ring~$R$ such that $H(G)=\Ll$~and~$H(BG)=\Sl$ are of the form
described in Section~\ref{Koszul-complex}.


\begin{proof}
  This computation can be found in~\cite[Example~2.2]{GugenheimMay:74}.
  (Their recursive formula looks slightly different because their
  {\cuponeproduct} is a right derivation.)
\end{proof}

We will choose specific representatives~$\xi_i'$ below.
The corresponding twisting cochain~$t$
then defines a functor $\rComod{C(BT)}\to\wComod\Sl$
which preserves homotopies:
If $f\colon M\to M'$ is a morphism of $C(BT)$-comodules, then
$f\otimes1\colon M\otimes_t\Ll\to M'\otimes_t\Ll$,
$m\otimes a\mapsto f(m)\otimes a$
is a morphism of right $\Ll$-modules, and analogously for homotopies.
Moreover, the $\Sl$-comodule structure on the homology of a $C(BT)$-comodule,
considered as weak $\Sl$-comodule, is the original $H(BT)$-comodule structure
because equation~\eqref{definition-module-structure-homology}
is how the cap product in~$H(M)$ by elements of~$H^*(BT)=\Su$ is defined.
We again suppress this functor from our notation, so that
$C(Y)\in\wComod\Sl$ for~$Y\in\Spaceover{BT}$.

We now explain how to choose the representatives~$\xi'_i$ we will work with.
In the case of a circle~$S^1$, we let $\xi'\in C^2(BS^1)$ correspond
under transgression to a cocycle ``dual'' to~$x'$:
It follows from the naturality of the cohomological
Leray--Serre spectral sequence~\eqref{Leray-Serre-cohomological}
that there exist
a cochain~$\xi'\in C^2(BS^1)$ and a cochain~$\chi\in C^1(ES^1)$
such that
\begin{subequations}\label{properties-chi}
\begin{align}
  \pair{x',i^*\chi} &= 1, \label{relation-x-chi}\\
  d\chi &= p^*\xi', \label{differential-chi}
\end{align}
\end{subequations}
where $i\colon S^1\to ES^1$ denotes the canonical inclusion over~$e_0$
and $p\colon ES^1\to BS^1$ the projection. Then $\xi'$ is a cocycle, and its
homology class~$\xi$ generates $H^*(BS^1)$.

For the torus~$T$ we define $\chi_i$ as the
cross product~$1\times\cdots\times\chi\times\cdots\times1$
corresponding to the $i$-th~factor of~$ET$, and analogously
for $\xi'_i$~and~$\xi_i$. Then the~$\xi_i$'s generate $\Su=H^*(BT)$ and are
represented by the~$\xi_i'$'s.

\subsection{The first natural transformation}\label{important-maps-1}

We now construct a map~$f\colon\K\to C(ET)$
which will help us to compare the algebraic and simplicial Koszul functors.
Here $\K=\K(P)$ is the Koszul complex as defined
in Section~\ref{Koszul-complex}.
For $r=1$~and~$l\in\N$ we recursively set
\begin{subequations}
\begin{align}
  f(1\otimes 1) &= e_0, \\
  f(x^l\otimes x) &= f(x^l\otimes1)\cdot x=f(x^l\otimes1)\cdot x', \\
  f(x^{l+1}\otimes1) &= S f(x^l\otimes x),
\end{align}
\end{subequations}
where $x\in\Sl$ denotes the canonical ``cogenerator'' dual to~$\xi$,
and $x'\in C(T)$ is the representative simplex of~$x\in\Ll$
chosen in the previous section.
For arbitrary~$r$, we compose this construction with the shuffle map,
\begin{equation}
  f\colon\K = \K(x_1)\otimes\dotsb\otimes \K(x_r)
    \to C(ES^1)\otimes\dotsb\otimes C(ES^1)\to C(ET),
\end{equation}
using equation~\eqref{BG-EG-products}.

\begin{proposition}\label{properties-f}
  This~$f$ is a morphism of $\Ll$-coalgebras
  and a strict morphism of weak $\Sl$-comodules.
\end{proposition}

The coalgebra structure on~$\K$ was defined in Section~\ref{Koszul-complex},
and $C(ET)$ is a $\Ll$-coalgebra
by \theoremref{C-X-C-G-coalgebra} and
\theoremref{section-Hopf-algebras}.
Recall that the second assertion of the
proposition\notheoremref[proposition]{properties-f} means that
\[
  f\otimes1\colon\K\otimes_P\Ll\to C(ET)\otimes_t\Ll,
  \qquad
  (s\otimes a)\otimes a'\mapsto f(s\otimes a)\otimes a'
\]
is a $\Ll$-equivariant chain map, where $t$ denotes the twisting cochain
constructed in the preceding section.
Here $\Ll$ acts of course only on the second factor of both tensor products.

\begin{proof}
  We may assume $r=1$ for the first claim because the shuffle map is a map of
  coalgebras and in addition equivariant if all but one space of a
  product have trivial group action, see \theoremref{shuffle-equivariant}.

  We start by verifying that $f$ is a chain map. By induction and
  equation~\eqref{S-homotopy}, we have for all~$l\in\N$
  (with the convention~$x^{-1}\otimes x:=0$)
  \begin{align*}
    d f(x^l\otimes1) &= d S f(x^{l-1}\otimes x)
      = f(x^{l-1}\otimes x)-S d f(x^{l-1}\otimes x) \\
    &= f(x^{l-1}\otimes x)-S f(d(x^{l-1}\otimes x))
      = f(x^{l-1}\otimes x), \\
    d f(x^l\otimes x) &= d(f(x^l\otimes1)\cdot x)
      =d f(x^l\otimes1)\cdot x=f(d(x^l\otimes1))\cdot x \\
    &= f(x^{l-1}\otimes x)\cdot x=f(x^{l-1}\otimes1)\cdot(x\wedge x)=0,
  \end{align*}
  which proves the claim by comparison with the explicit
  form~\eqref{Koszul-differential} of the differential on~$\K$.
  The $\Ll$-equivariance of~$f$ follows directly from the definition.
  To show that $f$ commutes with comultiplication, we have to establish the
  identities
  \begin{subequations}
  \begin{align}
    \AW\Delta_* f(x^l\otimes1)
      &= \sum_{m+n=l}f(x^m\otimes1)\otimes f(x^n\otimes1), \\
    \label{AW-Delta-f-x}
    \AW\Delta_* f(x^l\otimes x)
      &= \sum_{m+n=l}f(x^m\otimes x)\otimes f(x^n\otimes1) \\
      &+ \sum_{m+n=l}f(x^m\otimes1)\otimes f(x^n\otimes x). \notag
  \end{align}
  \end{subequations}
  They follow again by induction, \theoremref{cone-homotopy} and the fact
  that the Alexander--Whitney map is $\Ll$-equivariant. Hence $f$ is a
  map of left $\Ll$-coalgebras.

  That $f$ be a strict morphism of weak $\Sl$-comodules
  translates into the conditions 
  \begin{equation}\label{condition-morphism-comodules}
    \xi'_\mu\cap f(s\otimes a)=p^*\xi'_\mu\cap f(s\otimes a)
    =\begin{cases}
      f(\xi_i\cap s\otimes a) & \text{if~$\mu=\{i\}$}, \\
      0 & \text{otherwise}
    \end{cases}
  \end{equation}
  for all~$s\otimes a\in\K$ and~$\emptyset\ne\mu\subset[r]$,
  see equation~\eqref{morphism-weak-comodules}.
  Since we already know $f$ to be a map of coalgebras, the first alternative
  simplifies for~$r=1$ to~$\pair{p^*\xi',f(x\otimes1)}=1$.
  This holds by~\eqref{S-homotopy} and
  our choice~\eqref{properties-chi} of $\xi'$~and~$\chi$, because
  \begin{align*}
    \pair{p^*\xi',f(x\otimes1)} &= \pair{d\chi,S(e_0\cdot x')}
       =\pair{\chi,d S(e_0\cdot x')}\\
      &= \pair{\chi,e_0\cdot x'-S d(e_0\cdot x')}=\pair{i^*\chi,x'}=1.
  \end{align*}
  The case~$r>1$ now follows from the general identity
  \[
  (\alpha\times\beta)\cap\shuffle(a\otimes b)
    =\doublesign\beta a\shuffle(\alpha\cap a\otimes\beta\cap b)
  \]
  for chains $a$,~$b$ and cochains~$\alpha$,~$\beta$ on spaces $X$~and~$Y$,
  respectively, which is a consequence of \theoremref{shuffle-morph-coalgebras}
  and equation~\eqref{AW-shuffle-1}, \cf~\cite[\S 7.14]{Dold:80}.

  As for the second alternative of~\eqref{condition-morphism-comodules}, note
  that the \cuponeproduct s determining in \theoremref{functor-weak-S}
  the higher order elements~$\xi'_\mu$, $\degree\mu>1$, are here
  \crossoneproduct s of cochains coming from different factors of the induced
  decomposition of the classifying space~$BT$.
  Hence \theoremref{cross1-shuffle} gives the desired result.
\end{proof}

For a left $T$-space~$X$ we define the map~$\tcomp_X$ as
the bottom row of the following commutative diagram:
\begin{diagram*}
  & & \K\otimes C(X) & \rTo^{f\otimes1} & C(ET)\otimes C(X)
    & \rTo^\shuffle & C(ET\times X) \\
  & & \dTo & & \dTo & & \dTo>{{q_X}_*} \\
  \alg t C(X) & \rEqual & \K\otimesunder{\Ll}C(X)
    & \rTo & C(ET)\otimesunder{C(T)}C(X)
    & \rTo & C(ET\timesunder T X).
\end{diagram*}

\begin{corollary}\label{tcomp-X-multiplicative}
  This~$\tcomp_X$ is a morphism of coalgebras
  and a strict morphism of weak $\Sl$-comodules.
\end{corollary}

\begin{proof}
  The top row of the above diagram is a morphism of coalgebras
  by \theoremref{properties-f} and \theoremref{shuffle-morph-coalgebras},
  as are the vertical projections. This proves the first claim.
  The second follows from \theoremref{q-X-shuffle-equivariant}
  and \theoremref{properties-f}.
\end{proof}

We therefore obtain a natural transformation
\begin{equation}
  \newsymbol{\tcomp}\colon C\circ\Sim t\to\alg t\circ C.
\end{equation}
Note that $\tcomp_X$ may be written in the form
\begin{equation}
  \tcomp_X={q_X}_*\shuffle_{ET,X}(\thom\otimes1)\colon
    \alg t C(X)=\Sl\otimes_P C(X)\to C(\Sim t X)
\end{equation}
with~$\thom(s)=f(s\otimes1)$.
This shows the symmetry between $\tcomp$~and
the natural transformation~$\hcomp$
to be defined in the following section.

\subsection{The second natural transformation}\label{important-maps-2}

In this section we want to introduce a natural transformation
\begin{equation}
  \newsymbol{\hcomp}\colon\alg h\circ C\to C\circ\Sim h.
\end{equation}
by defining morphisms of $\Ll$-modules, natural in~$Y\in\Spaceover{BT}$,
\[
  \hcomp_Y\colon C(\Sim h Y)\to C(Y)\otimes_t\Ll=\alg h C(Y)
\]
of the form
\begin{equation}
  \hcomp_Y=(1\otimes\hhom)\AW_{Y,ET}{j_Y}_*
\end{equation}
where $\hhom\in\Hom_0(C(ET),\Ll)$ and
$j_Y\colon Y\timesover{BT}ET\to Y\times ET$ denotes the canonical inclusion.
We look at~$C(ET)$ as a \emph{left} $C(BT)$-comodule whose structure map
is the composition of the canonical map~$ET\to BT\times ET$
and the Alexander--Whitney map.
This allows to consider the
twisted homomorphism complex~$\Hom^t(C(ET),\Ll)$
as defined in Section~\ref{twisting}.

\begin{lemma}
  Such a~$\hcomp_Y$ is a morphism of right $\Ll$-modules
  if $\hhom$ is so and if it is a cycle in~$\Hom^t_0(C(ET),\Ll)$.
\end{lemma}

\begin{proof}
  This follows directly from the definitions and the associativity of
  the Alexan\-der--Whitney map.
\end{proof}

Our construction of such a~$\hhom$ goes as follows:
Recall that in Section~\ref{important-maps-1}
we have defined cochains~$\chi_i\in C^*(ET)$, $i\in[r]$,
satisfying
equations~\eqref{properties-chi}.
For a subset $\pi\subset[r]$ with at least two elements, we now set
\begin{equation}\label{definition-chi}
  \chi_\pi= \sign\pi\,\chi_{\pi^+}\cupone p^*\xi'_{\pi'},
\end{equation}
where $\pi^+$ again denotes the maximum of~$\pi$ and $\pi'$ its complement, and
\begin{equation}\label{definition-zeta}
  \zeta_\emptyset=1,\qquad
  \newsymbol{\zeta_\pi}=\sum_{\smashsubstack{(\mu,\nu)\vdash\pi\\\pi^+\in\mu}}
    \sign{(\nu,\mu)}\,\chi_\mu\cup\zeta_\nu
\end{equation}
for non-empty~$\pi$.
Note that we have~$\sign{\chi_\pi}=\sign{\zeta_\pi}=\sign\pi$.
We finally define
\begin{equation}
  \newsymbol{\hhom}=\sum_{\pi\subset[r]}x_\pi\otimes\zeta_\pi
    \in\Hom^t_0(C(ET),\Ll),
\end{equation}
where we have used
the isomorphism of graded $R$-modules~\eqref{isomorphism-Hom-tensor}.

\begin{proposition}
  This $\hhom$ is a $\Ll$-equivariant cycle.
\end{proposition}

\begin{proof}
  The condition~$d^{\,t}\hhom=0$ boils down to the relations
  \begin{equation*}
    d\zeta_\pi = -\sum_{\smashsubstack{(\mu,\nu)\vdash\pi \\ \mu\ne\emptyset}}
      \sign\mu\sign{(\nu,\mu)}\, p^*\xi'_\mu\cup\zeta_\nu
  \end{equation*}
  for all~$\pi\subset[r]$,
  \cf~the twisting cochain condition~\eqref{equation-a-pi}.
  They can be verified inductively using the formula
  \begin{equation*}
    \sign\pi\,d \chi_\pi = -p^*\xi'_\pi
      +\negthickspace\sum_{\substack{(\mu,\nu)\vdash\pi\\
          \pi^+\in\mu,\,\nu\ne\emptyset}}
        \negthickspace\sign{(\nu,\mu)}\,\chi_\mu\cup p^*\xi'_\nu
      -\negthickspace\sum_{\substack{(\mu,\nu)\vdash\pi\\
          \pi^+\in\nu,\,\mu\ne\emptyset}}
        \negthickspace\sign\nu\sign{(\nu,\mu)}\,p^*\xi'_\mu\cup\chi_\nu
  \end{equation*}
  for non-empty~$\pi$, which is a consequence of
  equation~\eqref{equation-a-pi} and the definition of~$\chi_\pi$.
  (Details for this and the following computations can be found
  in~\cite[Appendix~11]{Franz:01},
  again modulo somewhat different conventions.)

  The $\Ll$-equivariance of~$\hhom$ is equivalent to the identities
  \begin{equation*}
    x_i\cdot\zeta_\pi=\begin{cases}
      \sign{(\pi\setminus i,i)}\,\zeta_{\pi\setminus i}
        & \text{if~$i\in\pi$},\\
      0 & \text{otherwise}
    \end{cases}
  \end{equation*}
  for all~$i\in[r]$ and all~$\pi\subset[r]$.
  (The right $\Ll$-module structure of~$C(ET)$ dualises to a left structure
  on~$C^*(ET)$ by~$\pair{a\cdot\gamma,c}=\pair{\gamma,c\cdot a}$.)
  These equations follow by induction from
  \begin{equation}\label{x-c-pi}
    x_i\cdot\chi_\pi=\begin{cases}
      1 & \text{if~$\pi=\{i\}$},\\
      0 & \text{otherwise}.
    \end{cases}
  \end{equation}
  Let us prove the preceding line:
  We clearly have $x_i\cdot\chi_j=0$ for~$i\ne j$ because in this case
  the~$\chi_j$ comes from a factor of~$ET$ on which $x_i$ acts trivially.
  Similarly, $x_i\cdot\chi_\pi$ vanishes for~$\degree\pi\ge2$:
  The {\cuponeproduct} in~\eqref{definition-chi} is in fact
  a {\crossoneproduct} with the second factor
  coming from a trivial $T$-space, namely some~$BS^1$. Hence
  \[
    x_i\cdot \chi_\pi
    =-\sign\pi\,x_i\cdot\chi_{\pi^+}\cupone p^*\xi'_{\pi'}=0
  \]
  by Propositions \theoremref[proposition]{ST-equivariant}~%
  and~\theoremref[proposition]{cross1-zero}.
  It remains to show $x_i\cdot \chi_i=x\cdot\chi=1$.
  Note that $x\cdot\chi$ is a cocycle because
  \[
    d(x\cdot\chi)=-x\cdot d\chi=-x\cdot p^*\xi'=-p^*(x\cdot\xi')=0.
  \]
  Since $ET$ is connected, it suffices therefore to evaluate~$x\cdot\chi$
  on a single vertex. By our choice of~$\chi$ we obtain
  \[
    \pair{x\cdot\chi,e_0} = \pair{\chi,e_0\cdot x}=\pair{i^*\chi,x'}=1,
  \]
  which finally proves~\eqref{x-c-pi}.
\end{proof}

\subsection{The main theorem}\label{comparing-functors}

We can now compare the simplicial Koszul functors
(defined in Section~\ref{simplicial-functors})
\[
  \Sim t\colon\Space T\to\Spaceover{BT}
  \and
  \Sim h\colon\Spaceover{BT}\to\Space T
\]
with their algebraic counterparts (Section~\ref{Koszul-algebraic})
\[
  \alg t\colon\wMod\Ll\to\wComod\Sl
  \and
  \alg h\colon\wComod\Sl\to\wMod\Ll
\]
by using the natural transformations
\[
  \tcomp\colon\alg t\circ C \to C\circ\Sim t
  \and
  \hcomp\colon C\circ\Sim h \to \alg h\circ C
\]
constructed in the two preceding sections.
(Recall from Section~\ref{tori} that we consider the chain functor
as a functor~$\Space T\to\lMod{\Ll}$
and as a functor~$\Spaceover{BT}\to\wComod{\Sl}$.)

\begin{theorem}\label{natural-transformations-c-equivalences}
  The natural transformations $\tcomp$~and~$\hcomp$ are quasi-equivalences.
\end{theorem}

\begin{proof}
  Let $X$ be a left $T$-space. We want to use the Leray--Serre theorem to show
  that $\tcomp_X$ induces an isomorphism in homology.
  If we filter $\alg t C(X)=\Sl\otimes_{u_P}C(X)$ by $\Sl$-degree
  and $C(\Sim t X)$ as described in Section~\ref{fibre-bundles}, 
  then $\tcomp_X$ is filtration-preserving
  because for~$c\in C_q(X)$ and~$l\in\N$ the ``base component''
  of~${q_X}_*\shuffle\bigl(\thom(x^l)\otimes c\bigr)$
  is $q$-fold degenerate by the definition of the shuffle map.
  Hence $\tcomp_X$ induces a morphism from the spectral sequence associated
  to~$\alg t C(X)$ to the Leray--Serre spectral sequence
  for~$ET\timesunder T X$.
  By considering the inclusion~$X\hookrightarrow ET\timesunder T X$
  over~$b_0$, one sees
  that $E^2(\tcomp_X)$ is an isomorphism on the column~$p=0$.
  The map $E^2(\tcomp_X)$ being a morphism of $\Su$-modules,
  this extends to all~$p$.
  Hence $H(\tcomp_X)$ is an isomorphism, too.

  The proof for~$\hcomp$ is similar. Using the projection~%
  $Y\timesover{BT}ET\to Y$,
  we conclude that we have an isomorphism on the $E^1$ level for~$q=0$,
  which by $\Ll$-equivariance must hold for all~$q$.
\end{proof}


\subsection{Naturality}\label{naturality}

Up to now, we have kept the group~$G$ (resp.~$T$) fixed.
We now discuss what happens if one allows $G$ to vary.

Denote the category of all left modules by~$\LMod$.
A morphism in~$\LMod$ between an $A$-module~$M$ and an $A'$-module~$M'$
is a pair~$(\varphi,f)$, where $f\colon M\to M'$ is a chain map, equivariant
with respect to the morphism of algebras~$\varphi\colon A\to A'$.
The category~$\RComod$ of all right comodules is defined analogously.

On the side of spaces, one has the category~$\LSpace$ of all left group actions
and equivariant maps and the category~$\RSpaceover$ of all spaces over base
spaces.

Homotopies in these categories are defined similarly.
For instance, a homotopy between
two morphisms $(\varphi,f)$~and~$(\varphi,f')$ in~$\LMod$
is a $\varphi$-equivariant chain homotopy~$h$ between $f$~and~$f'$.
(Note that we require $f$,~$f'$
and~$h$ to be equivariant with respect to the same morphism of algebras.)

The following observation is a straightforward generalisation of
Propositions \theoremref[proposition]{chain-functor-Space-2-Mod}~and~%
\theoremref[proposition]{chain-functor-Spaceover-2-Comod}:

\begin{proposition}
  The chain functor~$C$ can be considered as
  homotopy-preserving
  functor~$\LSpace\to\LMod$ or~$\RSpaceover\to\RComod$.
\end{proposition}

For the simplicial Koszul functors one easily verifies the following:

\begin{proposition}
  The simplicial Koszul functors are natural with respect to morphisms
  in~$\LSpace$ on the one hand and morphisms in~$\RSpaceover$
  over maps~$B\varphi\colon BG\to BG'$ induced
  by morphisms of groups~$\varphi\colon G\to G'$ on the other.
\end{proposition}

Let $P'$ be another finitely generated free graded $R$-module,
also satisfying the conditions stated at the beginning of
Section~\ref{Koszul-complex}. It gives rise to an exterior algebra~$\Llprime$
and an symmetric coalgebra~$\Sl'$.
Any morphism of graded $R$-modules~$P\to P'$ canonically determines
a morphism of algebras~$\Ll\to\Llprime$ and
a morphism of coalgebras~$\Sl\to\Sl'$.
Call morphisms arising this way adapted,
as well as morphisms between (co)modules over
exterior algebras (resp.~symmetric coalgebras)
whose (co)algebra parts are of this type. Finally
extend this notion to weak (co)modules, whose morphisms
were defined in Sections \ref{weak-comodules}~and~\ref{weak-modules} 
as the morphisms of the Koszul dual strict (co)modules.
Then an adapted morphism of modules, considered as a morphism of weak
comodules, is still adapted, and similarly for adapted morphisms of comodules.

\begin{proposition}
  The algebraic Koszul functors are natural with respect to adapted morphisms
  of weak (co)modules.
\end{proposition}

This generality breaks down if one tries to compare
the simplicial Koszul functors with the algebraic ones.
The problem is that the definition of the
sweep action
as well as that of the twisting cochain~$t$ in Section~\ref{tori}
involve a decomposition of the torus~$T$ into circles
and furthermore the choice of representatives of some (co)homology classes.
A look at the construction of the connecting natural transformations
$\tcomp$~and~$\hcomp$
shows that we still have naturality in the following situations:

\begin{proposition} $ $\label{naturality-connecting-transformations}
  Let $\varphi\colon T\cong(S^1)^r\to T'\cong(S^1)^{r'}$ be a map of tori.
  Assume that it is componentwise
  with respect to the chosen decompositions in the sense that it just drops
  some factors and permutes the remaining ones, with the possible insertion
  of~$1$'s.
  \begin{enumerate}
  \item Let $X$ be a $T$-space, $X'$ a $T'$-space
    and $f\colon X\to X'$ a $\varphi$-equivariant map.
    Then, using the obvious notation, the following diagram commutes:
    \begin{diagram*}
      \alg t_PC(X) & \rTo^{\alg t f_*} & \alg t_{P'}C(X')\\
      \dTo<{\tcomp_X} & & \dTo>{\tcomp_{X'}}\\
      C(\Sim t_T X) & \rTo^{(\Sim t f)_*} & C(\Sim t_{T'}X').
    \end{diagram*}
  \item\label{naturality-monotone-maps}
    Let $Y$ be a space over~$BT$, $Y'$ a space over~$BT'$
    and $g\colon Y\to Y'$ a map over~$B\varphi\colon BT\to BT'$.
    If in addition $\varphi$ preserves the relative order
    of the surviving elements, then the following diagram commutes, too:
    \begin{diagram*}
      C(\Sim h_T Y) & \rTo^{(\Sim h g)_*} & C(\Sim h_{T'}Y')\\
      \dTo<{\hcomp_Y} & & \dTo>{\hcomp_{Y'}}\\
      \alg h_P C(Y) & \rTo^{\alg h g_*} & \alg h_{P'}C(Y').
    \end{diagram*}
  \end{enumerate}
  Here all horizontal arrows represent morphisms
  adapted to the morphism~$H_1(\varphi)\colon H_1(T)=P\to H_1(T')=P'$.
\end{proposition}

To give an example, the map~$(S^1)^4\to(S^1)^3$,
$(g_1,g_2,g_3,g_4)\mapsto (1,g_4,g_1)$ is componentwise, but not monotone.
The special cases $1\to T$~and~$T\to 1$ of the
proposition\notheoremref[proposition]{naturality-connecting-transformations}
were implicitly used in the proof
of~\theoremref{natural-transformations-c-equivalences}.


\section{Cohomology}\label{cohomology}

In this chapter we translate the two parts of our main
theorem~\theoremref[theorem]{natural-transformations-c-equivalences}
to cohomology.

\subsection{The singular Cartan model}\label{section-Cartan}

Let $T$ be a torus and $X$~a left $T$-space.
The left $\Ll$-action on~$C(X)$ introduced in Section~\ref{tori}
dualises to a right action on cochains
by~$\pair{\gamma,x_i\cdot c}=\pair{\gamma\cdot x_i,c}$. Let us replace this
by the left action
\begin{equation}\label{definition-left-action-cohomology}
  \pair{x_i\cdot\gamma,c}
  =-\sign\gamma\pair{\gamma\cdot x_i,c}
  =-\sign\gamma\pair{\gamma,x_i\cdot c}.
\end{equation}
(This is equation~\eqref{comparison-C-G-module-structures}
after taking homology.)

We write the complex dual
to~$\alg t C(X)=\Sl\otimes_P C(X)$ as~$\Su\otimes C^*(X)$
so that the differential takes the form
\begin{equation}\label{Cartan-differential}
  d(\sigma\otimes\gamma) =\sigma\otimes d\gamma
      +\sum_{i=1}^r\,\xi_i\sigma\otimes x_i\cdot\gamma.
\end{equation}
It follows from the properties of~$\Sl\otimes_P C(X)$ that its dual
is an $\Su$-module with $\Su$-bilinear product
\begin{equation}\label{Cartan-product}
 (\sigma'\otimes\gamma')(\sigma\otimes\gamma)
   = \sigma'\sigma\otimes\gamma'\cup\gamma.
\end{equation}
  
\theoremref{tcomp-X-multiplicative}
and \theoremref{natural-transformations-c-equivalences} now imply:

\begin{theorem}\label{singular-Cartan-model}
  Given the definitions \eqref{Cartan-differential}~and~\eqref{Cartan-product},
  the map
  \[
    \tcomp_X^*\colon C^*(ET\timesunder T X)\to\Su\otimes C^*(X)
  \]
  is a quasi-isomorphism of algebras.
  Moreover, the induced isomorphism in cohomology is $\Su$-equivariant.

  The map~$ \tcomp_X^*$ is natural with respect to morphisms between $T$-spaces
  and, more generally, with respect to maps~$X\to X'$ of spaces
  with torus actions which are equivariant with respect to
  componentwise morphisms~$T\to T'$.
\end{theorem}

Componentwise morphisms were explained
in \theoremref{naturality-connecting-transformations}.
Note that the correctness of the action of~$\Su=H_T^*(\oneptspace)$
already follows from the algebra structure and naturality.
In order to put the duals of weak $\Sl$-comodules into a categorical framework,
one could define the category of ``weak $\Su$-modules''.
The map~$\tcomp_X^*$ would then be a morphism in this new category.

That the singular Cartan model computes the equivariant cohomology algebra
of a $T$-space~$X$ was for~$T=S^1$ proven by
Jones~\cite[\S 3]{Jones:87} and Hood--Jones%
~\cite[Sec.~4]{HoodJones:87} in the context of cyclic homology.
For arbitrary~$T$, subanalytic $T$-spaces and real coefficients
the result (but without the algebra structure)
is due to Goresky, Kottwitz and MacPherson~%
\cite[Thm.~12.3]{GoreskyKottwitzMacPherson:98}.

Let us recall from Sections \ref{groups}~and~\ref{tori}
how the ``sweep action'' of~$\Ll$ on cochains comes about:
Let~$x_i\in H_1(T)$ be one of basis elements
corresponding to the chosen decomposition of~$T$ into circles.
It represents a loop~$x_i'$ at~$1\in T$ that winds once around
the $i$-th factor of~$T=(S^1)^r$.
Now let $c$ be an $n$-chain and $\gamma$ a cochain in~$X$.
Move $c$ around in~$X$ by applying the loop~$x_i'$ and triangulate
the result into  $(n+1)$-simplices as done
by the shuffle map~\eqref{shuffle}.
Evaluating $\gamma$ on this $(n+1)$-chain gives, up to the sign~$(-1)^n$,
the value of~$x_i\cdot\gamma$ on~$c$.

For $X$~a $T$-manifold and real coefficients one can deduce
the classical Cartan model from \theoremref{singular-Cartan-model}:
It is well-known that the restriction of cochains
\begin{equation}\label{restriction-smooth-simplices}
  C^*(X)\to C^*_\infty(X)
\end{equation}
to smooth singular simplices is a quasi-isomorphism of complexes,
as are the maps
\begin{equation}\label{T-invariant-form-to-cochain}
  \Omega(X)^T\hookrightarrow\Omega(X)\to C^*_\infty(X),
\end{equation}
where $\Omega(X)$ denotes the de~Rham complex of~$X$ and $\Omega(X)^T$
its $T$-invariants.
Letting $\Ll$ act on~$C^*_\infty(X)$ by the sweep makes
\eqref{restriction-smooth-simplices} a morphism of $\Ll$-modules.
On $T$-invariant forms $\Ll$ acts by contraction with
generating vector fields. It is not difficult to verify
that the composition~\eqref{T-invariant-form-to-cochain} is equivariant.
The point is that integrating a $T$-invariant differential form~$\gamma$
over~$x_i'\cdot c$ is the same as integrating $\gamma$, contracted by the
fundamental vector field generated by~$x_i'$, over~$c$.
The usual spectral sequence argument finally shows
that the classical Cartan model~$\Su\otimes\Omega(X)^T$
(with the same differential as above) is quasi-isomorphic as $\Su$-module
to its ``singular'' counterpart;
essentially the inverse argument was used
in~\cite[Sec.~18]{GoreskyKottwitzMacPherson:98}.
(For this to work one actually has to adjust some sign conventions
for differential forms.)
The direct comparison of the product structures
would of course be more involved. 

\subsection{The pull back of universal bundles}

Let $Y$ be a space over~$BT$. We write the complex dual
to~$\alg h C(Y)=C(Y)\otimes_t\Ll$ as~$C^*(Y)\otimes\Llstar$.
Then the differential is
\begin{equation}\label{gugenheim-may-differential}
  d(\gamma\otimes\alpha)=d\gamma\otimes\alpha
    -\sum_{\smashsubstack{\pi\subset\{1,\ldots,r\}\\\pi\ne\emptyset}}
      (-1)^{\degree\gamma+\degree\pi}\,
      \gamma\cup p^*(\xi_\pi')\otimes x_\pi\cdot\alpha
\end{equation}
(this is equation~\eqref{gugenheim-may-differential-introduction}
from the introduction),
and the left $\Ll$-action
\begin{equation}
  x_\pi\cdot(\gamma\otimes\alpha)
    =\doublesign\pi\gamma\,\gamma\otimes x_\pi\cdot\alpha.
\end{equation}
According to~\eqref{definition-left-action-cohomology}, $\Ll$ acts
on~$\Llstar$ by
\[
  x_j\cdot(\xi_{i_1}\wedge\cdots\wedge\xi_{i_k})
  =\xi_{i_2}\wedge\cdots\wedge\xi_{i_k}
  \qquad\hbox{if~$j=i_1$}.
\]

\begin{proposition}
  With these conventions,
  the map
  \[
    \hcomp_Y^*\colon C^*(Y)\otimes\Llstar\to C^*(Y\timesover{BT}ET)
  \]
  is a quasi-isomorphism of $\Ll$-modules.
  It is natural with respect to morphisms between $T$-spaces and,
  more generally, with respect to the morphisms described
  in \theoremref{naturality-monotone-maps}.
\end{proposition}

This sharpens the pertinent result
of Gugenheim--May~\cite[Ex.~2.2~\&~Thm.~3.3]{GugenheimMay:74}
by additionally describing the $\Ll$-action
(and also by allowing completely arbitrary~$R$).

\subsection{A quasi-isomorphism~$C^*(BT)\to H^*(BT)$}
\label{quasi-isomorphism-S-C-BT}

Taking~$X=\oneptspace$ a point in the singular Cartan model gives
together with the dual of the connecting morphism~$\tcomp_X$
from Section~\ref{important-maps-1} a quasi-isomorphism of algebras
\[
  \Bf^*\colon C^*(BT)\to\Su=H^*(BT).
\]
It turns out that $\Bf^*$ has another useful property:

\begin{proposition}\label{splitting-cup-one-0}
  The map~$\Bf^*$ is a quasi-isomorphism of algebras.
  Moreover, it annihilates all \cuponeproduct s.
\end{proposition}

The existence of such a map is also due to
Gugenheim and May~\cite[Thm.~4.1]{GugenheimMay:74}. Though somewhat
technical in nature, this result is of great importance, for instance to the
study of the cohomology of homogeneous spaces of Lie groups,
\cf~%
\cite{GugenheimMay:74}
or~\cite[\S 8.1]{McCleary:01} for instance.
The present construction of such a map is considerably
simpler than the original one given in the appendix to~\cite{GugenheimMay:74}.
Before trying to prove this result, I have checked some examples
with the help of the ``Kenzo'' program~\cite{DoussonSergeraert:98}.

\begin{proof}
  We actually need not invoke
  \theoremref{natural-transformations-c-equivalences}
  (as implicitly done above) to see that $\Bf^*$ induces
  an isomorphism in homology: This simply follows from the fact that we have
  chosen the~$\xi'_i$ as representatives of the~$\xi_i\in\Su=H^*(BT)$ in
  Section~\ref{tori}.

  That $\Bf^*$ annihilates all \cuponeproduct s is equivalent to the
  vanishing of
  \begin{align*}
    \ST_{BT,BT}{\Delta_{BT}}_*\Bf &\colon \Sl\to C(BT)\otimes C(BT) \\
  \intertext{and to that of}
    (p_*\otimes p_*)\ST_{ET,ET}{\Delta_{ET}}_*f
    &\colon\K\to C(BT)\otimes C(BT),
  \end{align*}
  where $f\colon\K\to C(ET)$ is the map constructed
  in Section~\ref{important-maps-1}.
  We actually prove the stronger statement that
  \[
    \ST_{BT,ET}\Deltapp_*f
     = (p_*\otimes1)\ST_{ET,ET}{\Delta_{ET}}_*f
    \colon\K\to C(BT)\otimes C(ET)
  \]
  vanishes, where $\Deltapp$ is the canonical map~$ET\to BT\times ET$.
  We proceed by double induction on the rank~$r$ of~$T$ and the degree
  of~$c=x^\alpha\otimes x_\pi\in\K$, the case~$r=0$ being trivial.
  If $r>0$~and~$\pi$ non-empty, then $x_\pi=a\wedge x_i$ for some~$a\in\Ll$
  and some~$i$. By \theoremref{ST-equivariant-1},
  the Steenrod map~$\ST_{BT,ET}$ is equivariant with respect to multiplication
  by~$x_i$ because the latter is of degree~$1$. Hence
  \begin{align*}
    \ST\Deltapp_*f(c)
    &= \ST\Deltapp_*\bigl(f(x^\alpha\otimes a)\cdot x_i\bigr)
     = \ST\bigl(\Deltapp_*f(x^\alpha\otimes a)\cdot x_i\bigr) \\
    &= \ST\Deltapp_*f(x^\alpha\otimes a)\cdot x_i
     = 0
  \end{align*}
  by induction.

  It remains the case~$\pi=\emptyset$,
  \ie, $c=x^\alpha\otimes1$. We may assume all~$\alpha_i>0$.
  (Otherwise use the result for smaller~$r$).
  Formula~\eqref{shuffle-S} shows that the cross product of two
  chains lying in the image of the respective cone operators does so itself.
  This generalises readily to several factors and applies therefore to
  \[
    f(c)=\shuffle_{ES^1,\ldots,ES^1}\bigl(
           f(x_1^{\alpha_1}\otimes1)\otimes\cdots
             \otimes f(x_r^{\alpha_r}\otimes1)\bigr).
  \]
  Since $S$ is a contracting homotopy by equation~\eqref{S-homotopy}
  and $SS=0$, we conclude that
  \[
    f(c) = {S d}f(c)+d S f(c) = Sf(dc).
  \]
  Applying equation~\eqref{ST-S} yields
  \[
    \ST\Delta_*f(c)
    = \ST\Delta_*S f(dc)
    = (S\otimes S)\AW\Delta_*f(dc)-(1\otimes S)\ST\Delta_*f(dc).
  \]
  The first summand vanishes because $f$ is morphism of coalgebras and~$SS=0$,
  \cf~the explicit form~\eqref{AW-Delta-f-x} of~$\AW\Delta_*f$.
  We project the remaining term to~$C(BT)\otimes C(ET)$:
  \begin{align*}
    \ST\Deltapp_*f(c)
    &= -(p_*\otimes1)\ST\Delta_*f(c)
     = -(1\otimes S)(p_*\otimes1)\ST\Delta_*f(d c) \\
    &= -(1\otimes S)\ST\Deltapp_*f(d c) = 0,
  \end{align*}
  again by induction. This finishes the proof.
\end{proof}

\section{Intersection homology}\label{intersection-homology}

\subsection{Motivation}\label{motivation-allowable}

This section serves as a motivation of our definition of
allowable subsets of a simplicial set in the next section.

Let $X$ be a stratified pseudomanifold and denote by~$X_k$ the union of
its strata of codimension at least~$k$.
Let $\T$ be a triangulation of~$X$ compatible with the stratification.
Assume furthermore that the vertices of~$\T$ are partially ordered
in such a way that the vertices of each simplex~$\sigma\in\T$ are ordered and
the intersection of~$\sigma$ with any~$X_k$ is a \newterm{back face}
of~$\sigma$, \ie, spanned by the last $\dim(\sigma\cap X_k)+1$~vertices
of~$\sigma$ in the given ordering.
Guided by~\cite{GoreskyMacPherson:86}, we call such
a vertex-ordered triangulation \newterm{flaglike}.
For example, the barycentric subdivision of
any triangulation compatible with the stratification is flaglike
if one orders the vertices by decreasing dimension of the faces of which
they are the barycentres.

Let $\p=(\p_0=\p_1=\p_2=0,\p_3,\ldots)$ be a perversity.
Call a simplex~$\sigma\in\T$ allowable if for all~$k$
the intersection~$\sigma\cap X_k$ has codimension at least~$k-\p_k$ in~$\sigma$,
and a chain in~$C(\T)$ allowable if it is a linear combination
of allowable simplices.
(We just say ``allowable'' because we will not deal with more than one
perversity at a time.)
Then define $I^{\p}C(\T)\subset C(\T)$ as the subcomplex of all allowable
chains with allowable boundaries.
Goresky and MacPherson~\cite{GoreskyMacPherson:86} have shown that
provided $\T$ is flaglike (as we assume here), the homology 
of~$I^{\p}C(\T)$ is the usual
intersection homology~$I^{\pneg}H(X)$ (with compact support)
as defined in~\cite{GoreskyMacPherson:80}.

The complex~$I^{\p}C(\T)$ has the nice property that it can easily
be endowed with the structure of a right $C(\T)$-comodule:
Since any intersection~$\sigma\cap X_k$ is a back face
of~$\sigma\in\T$, all \newterm{front faces}~$\tpartial^i\sigma$
of an allowable~$\sigma$ are again allowable. In other words, the set of
allowable simplices is stable under the last face map of the corresponding
simplicial set. Consequently, the image of a chain~$c\in I^{\p}C(\T)$ under the
Alexander--Whitney map~\eqref{AW} is of the form
\begin{equation}\label{image-AW-allowable}
  \sum_i c'_i\otimes\sigma_i''
\end{equation}
for some non-degenerate simplices~$\sigma_i''$
and some allowable chains~$c'_i$.
The Alexander--Whitney being a chain map, this already implies that it maps
$I^{\p}C(\T)$ to~$I^{\p}C(\T)\otimes C(\T)$:
Since $C(\T)$ is free over~$R$, the chain~$\AW(c)$ lies
in~$I^{\p}C(\T)\otimes C(\T)$ if and only if $(d\otimes1)\AW(c)$ is also
of the form~\eqref{image-AW-allowable}. This is true because
\[
  (d\otimes1)\AW(c)=d\AW(c)-(1\otimes d)\AW(c)=\AW(d c)-(1\otimes d)\AW(c).
\]
Hence $I^{\p}C(\T)$ is a right $C(\T)$-comodule,
and $I^{\pneg}H(X)$ a right $H(X)$-comodule.

For compact~$X$, the induced dual $H^*(X)$-module structure
on~$I^{\pneg}H(X)$ is the usual one, which is defined as the
intersection product~$\pi([\gamma])\times[c]$
of the Poincar\'e dual~$\pi([\gamma])\in I^{\pzero}H(X)$
of~$[\gamma]\in H^*(X)$
with~$[c]\in I^{\pneg}H(X)$.
This follows from inspection of the definition of the intersection product
of two transversal chains~\cite[\S 2.1]{GoreskyMacPherson:80} and the fact
that the Poincar\'e dual chain~$\pi(\gamma)=\gamma\cap[X]$
of~$\gamma\in C^*(\T)$, which is a union of simplices of the barycentric
subdivision of~$\T$, is transversal to all strata,
\cf~\cite[\S\,3.B]{Brasselet:95}.

\medbreak

Of course, if some Lie group~$G$ acts on~$X$ compatibly
with the stratification, we also want the complex of intersection chains
to be a $C(G)$-module.
If $G$ is not finite, then a triangulation is not convenient.
For this reason we switch to singular chains and use a variant
of King's definition of intersection homology \cite{King:85}:

Call a singular simplex~$\sigma\colon\Delta_n\to X$ allowable
if the inverse image of~$X_k$ is contained in the $(k-\p_k)$-codimensional
back face of~$\sigma$ for all~$k$. The graded set~$V$ of allowable singular
simplices is again stable under the last face map.
Imitation of the above definition leads to the
subcomplex~$C(V\subset X)\subset C(X)$ of normalised allowable chains
with allowable boundaries. It is a $C(X)$-comodule for the same reason
as before. 
In addition, any degeneration~$s_i\sigma$ of an allowable
simplex~$\sigma$ is again allowable. (This of course is true for
a flaglike triangulation~$\T$ as well.)
Hence if some Lie group~$G$ acts on~$X$ compatibly with the stratification,
we can use the shuffle map to define a $C(G)$-module structure
on~$C(V\subset X)$ as in Section~\ref{groups}.

Let $Y$ be a manifold and give $X\times Y$ the stratification
induced by that of~$X$.
Since the Eilenberg--Mac\,Lane maps~\eqref{eilenberg-zilber-maps}
for the pair~$(X,Y)$ only apply
$\tpartial$~and degeneracy maps to simplices from~$X$, they restrict to maps
between $C(V\subset X)\otimes C(Y)$~and~$C(V\times Y\subset X\times Y)$,
and properties~\eqref{Eilenberg-Zilber-relations} still hold.
We therefore get a K\"unneth theorem for~$H(V\times Y\subset X\times Y)$.%
\footnote{We will see in a minute that $H(V\times Y\subset X\times Y)$
is the usual intersection homology of~$X\times Y$.
We therefore prove the conjecture made in~\cite[p.~152]{King:85}.}
If $\U$ is a triangulation of~$Y$, then the same holds for
$I^{\p}C(\T)\otimes C(\U)$~and~$I^{\p}C(\T\times\U)$,
where $\T\times\U$ denotes the vertex-ordered triangulation
of~$X\times Y$ obtained by triangulating
each~$\sigma\times\tau$, $\sigma\in\T$,~$\tau\in\U$ as done by the shuffle map.
Note that $\T\times\U$ is again flaglike.

Since the vertices of any simplex in a flaglike triangulation~$\T$ of~$X$
are ordered, we have a canonical
chain map~$\alpha_\T\colon I^{\p}C(\T)\to C(V\subset X)$.

\begin{lemma}
  The induced map~$I^{\pneg}H(X)\to H(V\subset X)$ is an isomorphism
  of $H(X)$-comodules
  and does not depend on the flaglike triangulation of~$X$.
\end{lemma}

\begin{proofnoqed}
  The map~$\alpha_\T$ is clearly a morphism of comodules,
  equivariant with respect to the quasi-isomorphism
  of coalgebras~$C(\T)\to C(X)$.
  Hence the induced map in homology is a morphism of $H(X)$-comodules.

  Let us show the independence of the triangulation next:
  Let $\T$ and $\T'$ be two flaglike triangulations of~$X$.
  Since they possess a common flaglike refinement
  (e.g., the barycentric subdivision of any refinement), we may assume
  $\T'$ to be a refinement of~$\T$.
  The same example shows that we may also assume that the supporting
  simplices~$\sigma_1$,~$\sigma_2\in\T$ of two vertices~$v_1$,~$v_2\in\T'$
  satisfy $\sigma_1\subset\sigma_2$ if~$v_1\le v_2$.

  Denote by~$i$ the canonical
  inclusion~$I^{\p}C(\T)\hookrightarrow I^{\p}C(\T')$ which sends each
  $n$-simplex~$\sigma\in\T$ to the sum of the $n$-simplices in~$\T'$
  contained in the support of~$\sigma$.
  Then $H(i)$ is an isomorphism by~\cite{GoreskyMacPherson:86}.
  It therefore suffices to construct for each cycle~$c\in I^{\p}C(\T)$
  a chain~$b\in C(V\subset X)$ such that~$d b=\alpha_\T(c)-\alpha_{\T'}(i(c))$.
  This can be done analogously to the construction of~$\beta$
  in~\cite{GoreskyMacPherson:86}. Here it is crucial that $C(V\subset X)$
  is a subcomplex of the \emph{normalised} chain complex of~$X$;
  one also needs the condition ``$v_1\le v_2\implies\sigma_1\subset\sigma_2$''
  mentioned above.

  In order to prove that the map~$I^{\pneg}H(\T)\to H(V\subset X)$
  is an isomorphism, we verify the conditions for
  King's comparison theorem~\cite[Thm.~10]{King:85}:
  \begin{enumerate}
  \item The usual proof of the exactness of the Mayer--Vietoris sequence
    by barycentric subdivision (with explicit homotopy
    as in~\cite[Appendix~1]{Vick:73}) works for~$H(V\subset X)$
    if the cone operator \emph{prepends} the new vertex to a simplex.
  \item For~$Y$,~$\U$,~$\T\times\U$~as above and $W$ the set of allowable
    simplices in~$X\times Y$, we get a commutative diagram
    \begin{diagram*}
      I^{\p}C(\T)\otimes C(\U) & \rTo^\shuffle & I^{\p}C(\T\times\U) \\
      \dTo<{\alpha_\T\otimes\alpha_\U} & & \dTo>{\alpha_{\T\times\U}} \\
      C(V\subset X)\otimes C(Y) & \rTo^\shuffle & C(W\subset X\times Y).
    \end{diagram*}
    Here both shuffle maps and~$\alpha_\U$ induce isomorphisms in homology,
    hence if $\alpha_\T$ does so, too, then also $\alpha_{\T\times\U}$.
    Since $\T\times\U$ is flaglike, the latter map
    is $I^{\pneg}H(X)\to H(V\subset X)$.
  \item The usual computation of the intersection homology of a cone
    carries over to~$H(V\subset X)$ if one adds the apex to a simplex
    as the \emph{last} vertex.
  \item is trivial.
  \item follows from the K\"unneth theorem discussed above.\qedhere{83pt}
  \end{enumerate}
\end{proofnoqed}

A payoff of the use of singular intersection chains is that this does not
require the space to be finite-dimensional -- not even for the $H^*(X)$-action.
The only ingredients we need
are a decreasing filtration of~$X$ by (possible countably many)
subsets~$X=X_0\supset X_1\supset X_2\supset\cdots$
and a ``generalised perversity''~$\p\colon\N\to\N$
with only $\p_0=0$~and~$\p_{j+1}\le\p_j+1$.
But in this generality one certainly loses the independence of the filtration
and the ``Poincar\'e'' duality between complementary intersection homology
groups, though one might hope to define a cohomological
``intersection cup product''.

If $X$ is a $G$-space with a $G$-stable filtration, we get an induced
filtration~%
$
  EG\timesunder G X=EG\timesunder G X_0\subset EG\timesunder G X_1\subset\cdots
$ of the Borel construction,
hence a set of allowable simplices~$V_G$ in~$EG\timesunder G X$,
and we can define $H^G(V\subset X)=H(V_G\subset EG\timesunder G X)$
without using finite-dimensional approximations 
to~$EG\to BG$.
We will do this in the simplicial setting in Section~\ref{main-allowable}.
For the moment being, we only note the identity
\[
  H^G(V\subset X)=I^\p H^G(X):=\lim_{n\to\infty}I^\p H(X\timesunder G EG_n), 
\]
where the~$EG_n$ is an $n$-universal approximation to~$EG$.
This is a consequence of the following more general observation:

\begin{proposition}
  Let $X$ be the direct limit of an increasing sequence of subspaces~$X_{(n)}$
  such that points are closed in~$X$. Suppose that $X$ has a filtration
  as above and give each~$X_{(n)}$ the induced filtration. Then for each
  generalised perversity~$\p$ one has
  \[
    H(V\subset X)=\lim_{n\to\infty}H(V_{(n)}\subset X_{(n)}).
  \]
\end{proposition}

\begin{proof}
  Since direct limits commute with homology, it suffices to show that
  $C(V\subset X)$ is the direct limit of the~$C(V_{(n)}\subset X_{(n)})$.
  This follows from the fact that, being compact, the image of a
  singular simplex~$\sigma\colon\Delta_k\to X$ already lies
  in some~$X_{(n)}$.
\end{proof}

\subsection{Allowable subsets}

We now abstract from the preceding discussion and define an
\newterm{allowable subset}~$V\subset X$ of a space (i.e., simplicial set)~$X$
to be a collection of subsets~$V_n\subset X_n$ which is stable
under all degeneracy maps and under the last face map~$\tpartial$.
If $X$ is a $G$-space, we also require all $V_n$ to be $G_n$-stable.
A simplex in~$X$ is called \newterm{$V$-allowable} if it lies in~$V$.
A chain in~$X$ is called \newterm{$V$-allowable} if it is a linear combination
of $V$-allowable simplices. We denote by $C(V\subset X)\subset C(X)$
the subcomplex of normalised $V$-allowable chains
with $V$-allowable boundaries.
Analogous to the previous section, $C(V\subset X)$ is a right $C(X)$-comodule,
and also a left $C(G)$-module if $X$ is a left $G$-space.

The definition of~$C(V\subset X;N)\subset C(X;N)$ with coefficients in
an $R$-module~$N$ is analogous. Note that
$C(V\subset X;N)\ne C(V\subset X)\otimes N$ in general,
but equality holds if $N$ is free over~$R$.

Let $W\subset Y$ be another allowable subset and $f\colon X\to Y$
a map of spaces. If $f$ maps $V$-allowable simplices to $W$-allowable ones,
then it induces a chain map~$f_*\colon C(V\subset X)\to C(W\subset Y)$.
The analogous statement holds for homotopies. We call such maps
\newterm{allowable}.


\subsection{The main theorem for allowable subsets}\label{main-allowable}

Let $B$ be a space, $\tau\colon B_{>0}\to G$ a twisting function
and $V\subset F$ an allowable subset of a left $G$-space~$F$.
Then $B\times_\tau V=(B_n\times V_n)_n\subset B\times_\tau F$ is allowable
because $V$ is $G$-stable.
Likewise, if $W\subset B$ is allowable, then so is
$W\times_\tau F\subset B\times_\tau F$.

This implies that the simplicial Koszul functors~$\Sim t$,~$\Sim h$ extend
to allowable subsets of $G$-spaces and spaces over~$BG$, respectively.
Moreover, the morphisms $\Sim P_X$~and~$\Sim I_Y$
defined by~\eqref{definitions-P-I} are allowable, hence restrict to maps
\begin{equation*}
  \Sim h\Sim t(V\subset X)\to V\subset X
  \and
  W\subset Y\to \Sim t\Sim h(W\subset Y),
\end{equation*}
which are quasi-isomorphisms because the homotopy inverses as well as
the homotopies given in the proof of \theoremref{simplicial-composition}
are allowable, too.

The natural transformations $\tcomp$~and~$\hcomp$
from Sections \ref{important-maps-1}~and~\ref{important-maps-2}
are also well-defined
in this new setting: The map~$\tcomp_X$ applies only degeneracy maps to
simplices of a $T$-space~$X$, and $\hcomp_Y$ only the last face map to
simplices of a space~$Y$ over~$BT$.

In order to prove the generalisation
of \theoremref{natural-transformations-c-equivalences} to allowable subsets,
we finally need an appropriate version of the Leray--Serre theorem:

\begin{proposition}\label{Leray-Serre-allowable-base}
  Let $B\times_\tau F$ be a fibre bundle and let $W\subset B$ be allowable.
  If $H(F)$ is free over~$R$,
  then the filtration of~$W\subset B$ by the skeletons of the base
  leads to a spectral sequence with
  \[
    E^2_{p q}(W\subset B, F)=H_p(W\subset B; H_q(F)).
  \]
  This is an isomorphism of $H(G)$-modules if~$F=G$.
\end{proposition}

The reason for assuming $H(F)$ to be free is to assure
that $C(W\subset B; H_q(F))$
equals~$C(W\subset B)\otimes H_q(F)$.

\begin{proof}
  The simplicial Leray--Serre theorem as described
  in Section~\ref{fibre-bundles} can be deduced from
  the twisted Eilenberg--Zilber theorem~\cite{Shih:62}.
  (The algebraic essence of it~\cite{Brown:65} is nowadays called the
  ``basic perturbation lemma''.)
  All this goes through in the present setting:

  As discussed in Section~\ref{motivation-allowable},
  all three Eilenberg--Mac\,Lane maps restrict to maps
  between $C(W\times F\subset B\times F)$~and~$C(W\subset B)\otimes C(F)$.
  The complexes $C(W\times F\subset B\times F)$~%
  and~$C(W\times_\tau F\subset B\times_\tau F)$
  have the same underlying graded $R$-modules
  because the allowable simplices are the same
  and the difference~$q=d_\tau-d$ of the two differentials
  maps allowable chains to allowable chains.
  This implies in particular that $q$ is a well-defined
  endomorphism of~$C(W\times F\subset B\times F)$.
  Application of the basic perturbation lemma now
  shows that $C(W\times_\tau F\subset B\times_\tau F)$
  is homotopy equivalent to~$C(W\subset B)\otimes_u C(F)$
  for the same twisting cochain~$u=u(\tau)$ as usual.
  This homotopy equivalence is not one of $C(B)$-comodules any more,
  but still one of $C(G)$-modules if~$F=G$.
  Filtering $C(W\times_\tau F\subset B\times_\tau F)$
  as in Section~\ref{fibre-bundles} and $C(W\subset B)\otimes_u C(F)$
  by the degree of the first factor gives the result, \cf~\cite[\S 32]{May:68}.
\end{proof}


\begin{theorem} $ $\label{allowable-subsets-theorem}
  \begin{enumerate}
  \item Let $X$ be a $T$-space and $V\subset X$~allowable.
    Then $\tcomp_X$ restricts to a
    quasi-isomorphism of weak $\Sl$-comodules
    \[
      \tcomp_X\colon \alg t C(V\subset X)
        \to C\bigl(\Sim t(V\subset X)\bigr).
    \]
  \item Let $Y$ be a space over~$BT$ and $W\subset Y$ allowable.
    Then $\hcomp_Y$ restricts to a quasi-isomorphism of $\Ll$-modules
    \[
      \hcomp_Y\colon C\bigl(\Sim h(W\subset Y)\bigr)
        \to \alg h C(W\subset Y).
    \]
  \end{enumerate}
\end{theorem}

\begin{proof}
  For~$\hcomp$ the proof of \theoremref{natural-transformations-c-equivalences}
  carries over without changes if one uses the new version of the Leray--Serre
  theorem proven above.

  In order to show that~$\tcomp_X$ is a quasi-isomorphism for all
  allowable~$V\subset X$, we proceed in a somewhat roundabout way:
  We note first that it is sufficient to prove the assertion for
  fibre products. This follows from the commutative diagram
  of weak $\Sl$-comodules
  \begin{diagram*}
    \alg t C\bigl(\Sim h\Sim t (V\subset X)\bigr)
      & \rTo & \alg t C(V\subset X)\\
    \dTo<{\tcomp_{\Sim h\Sim t X}} & & \dTo>{\tcomp_X}\\
    C\bigl(\Sim t\Sim h\Sim t (V\subset X)\bigr)
      & \rTo & C\bigl(\Sim t(V\subset X)\bigr),
  \end{diagram*}
  whose horizontal arrows are quasi-isomorphisms by the remarks above
  together with \theoremref{Koszul-preserve-quasi-isomorphism}.

  For~$V\subset X=\Sim h (W\subset Y)$ we use the naturality
  of our constructions with respect
  to the morphism of groups~$T\to1$. We write $\Sim h=\Sim h_T$ in order to
  distinguish it from the trivial functor~$\Sim h_1$, and similarly for
  all other functors. The commutative diagram of chain maps
  \begin{diagram*}
      C\bigl(\Sim t_T\Sim h_T(W\subset Y)\bigr)
        & \rTo & C\bigl(\Sim t_T\Sim h_1 (W\subset Y)\bigr)
        & \rTo & C\bigl(\Sim t_1\Sim h_1 (W\subset Y)\bigr) \\
      \uTo<{\tcomp_{T,\Sim h_T Y}} & & \uTo<{\tcomp_{T,\Sim h_1 Y}}
        & & \uTo<{\tcomp_{1,\Sim h_1 Y}} \\
      \alg t_T C\bigl(\Sim h_T (W\subset Y)\bigr)
        & \rTo & \alg t_T C\bigl(\Sim h_1 (W\subset Y)\bigr)
        & \rTo & \alg t_1 C\bigl(\Sim h_1 (W\subset Y)\bigr) \\
      \dTo<{\alg t_T\hcomp_{T,Y}} & & \dTo<{\alg t_T\hcomp_{1,Y}}
        & & \dTo<{\alg t_1\hcomp_{1,Y}} \\
      \alg t_T\alg h_T C(W\subset Y)
        & \rTo & \alg t_T\alg h_1 C(W\subset Y)
        & \rTo & \alg t_1\alg h_1 C(W\subset Y)
  \end{diagram*}
  condenses to
  \begin{diagram*}
      C\bigl(\Sim t\Sim h (W\subset Y)\bigr) \\
      \uTo<{\tcomp_{\Sim h Y}} & \rdTo \\
      \alg t C\bigl(\Sim h (W\subset Y)\bigr) & & C(W\subset Y). \\
      \dTo<{\alg t\hcomp_Y} & \ruTo \\
      \alg t\alg h C(W\subset Y)
  \end{diagram*}
  Since we know all maps but $\tcomp_{\Sim h Y}=\tcomp_X$ to induce
  isomorphisms in homology, this map must do so, too.
\end{proof}

\begin{corollary}
  Using the same notation as in the
  theorem\notheoremref[theorem]{allowable-subsets-theorem}
  and moreover the differentials
  \eqref{Cartan-differential}~and~\eqref{gugenheim-may-differential},
  we have:
  \begin{enumerate}
  \item The map
    \[
      H^*(\tcomp_X)\colon H^*\bigl(\Sim t(V\subset X)\bigr)
        \to \Su\otimes H^*(V\subset X)
    \]
    is an isomorphism of $\Su$-modules.
  \item The map
    \[
      H^*(\hcomp_Y)\colon H^*(W\subset Y)\otimes\Llstar
        \to H^*\bigl(\Sim h(W\subset Y)\bigr)
    \]
    is an isomorphism of $\Ll$-modules.
  \end{enumerate}
\end{corollary}

\bibliographystyle{plain}
\bibliography{abbrev,algtop}

\end{document}